# On the variational method for Euclidean quantum fields in infinite volume[*]

N. Barashkov

Department of Mathematics and Statistics
University of Helsinki

*Email:* nikolay.barashkov@helsinki.fi

M. Gubinelli

Mathematical Institute
University of Oxford

*Email:* gubinelli@maths.ox.ac.uk

**Abstract**

We investigate the infinite volume limit of the variational description of Euclidean quantum fields introduced in a previous work. Focussing on two dimensional theories for simplicity, we prove in details how to use the variational approach to obtain tightness of $\varphi_2^4$ without cutoffs and a corresponding large deviation principle for any infinite volume limit. Any infinite volume measure is described via a forward–backwards stochastic differential equation in weak form (wFBSDE). Similar considerations apply to more general $P(\varphi)_2$ theories. We consider also the $\exp(\beta\varphi)_2$ model for $\beta^2 < 8\pi$ (the so called full $L^1$ regime) and prove uniqueness of the infinite volume limit and a variational characterization of the unique infinite volume measure. The corresponding characterization for $P(\varphi)_2$ theories is lacking due to the difficulty of studying the stability of the wFBSDE against local perturbations.

**A.M.S. subject classification:** 60H30, 81T08

The aim of this paper is to develop further the variational approach to Euclidean quantum fields initiated in [3] by studying the infinite volume limit of the variational description and some of its consequences. For simplicity we will focus on the two dimensional setting. We will highlight some unsolved issues which would have to be worked out in this simpler setting before moving to more challenging models. The infinite volume limit in the $\varphi_3^4$ theory considered in [3] can be attacked along the lines suggested here.

The key tool is a novel kind of stochastic equations for Euclidean fields which involve an additional one dimensional parameter. This parameter already appeared in [3] linked to a continuous scale decomposition of the relevant random fields. Here the additional parameter $t$ can be chosen to take values in $[0, 1]$ with no particular association to the scale of fluctuations thanks to the milder irregularity of the two dimensional Gaussian free field and to the simpler renormalization (i.e. Wick ordering) needed. In a general perspective one can interpret the variational method and the associated stochastic equations as a way to *stochastically quantize* a given field theory, alternative to the classical approach of Parisi–Wu via a Langevin equation [32] and to the elliptic stochastic quantization recently introduced in [2, 1] via the supersymmetric dimensional reduction of Parisi–Sourlas [31].

The variational method we introduced in [3] has been recently used in various contexts to provide novel results: Bringmann used it to construct singular Gibbs measures [11] and study their invariance under dispersive dynamics [12], Chandra, Gunaratnam and Weber used it to show the existence of phase transition in the $\varphi_3^4$ model [13], Oh, Okamoto and Tolomeo [30] to prove existence and non-existence of certain Gibbs measures. In [5] we provide another approach to describe the most singular part of the $\varphi_3^4$ measure by approximating the optimal drift of the variational description. The PhD thesis of the first author [4] describes also how to study the infinite volume limit of the Sine–Gordon model via the variational method including a full verification of the Osterwalder–Schrader axioms.

---

[*]. This article has been written using GNU T<sub>E</sub>X<sub>MACS</sub> [22].

Similar ideas relying on a coupling between the Gaussian free field and an interacting field have been used by Bauerschmidt and Hofstetter to control the maximum of the Sine–Gordon model [6]. More generally, mapping a Brownian motion into a measure of interest has revealed itself an useful method also in the context of functional inequalities, see e.g. the work of Borell [8], Lehec [28] and the recent paper of Mikulincer and Shenfeld [29] on the "Brownian transport map".

Let us now give an informal and heuristic picture of the results of this paper. The basic object of interest is a two dimensional Euclidean quantum field [34, 17], namely a probability measure $\theta$ on the Schwarz space $\mathcal{S}'(\mathbb{R}^2)$ of tempered distributions given by the expression

$$\theta(d\varphi) = \frac{e^{-V(\varphi)}}{Z}\theta^0(d\varphi),$$

where $\theta^0$ is the Gaussian free field (GFF) on $\mathbb{R}^2$ with covariance $(m^2 - \Delta)^{-1}$ and $V$ is a local interaction term given by

$$V^\xi(\varphi) = \int_{\mathbb{R}^2} \xi(x) p(\varphi(x)) dx,$$

with $\xi\colon \mathbb{R}^2 \to \mathbb{R}$ a compactly supported function to make the spatial integral finite and $p$ is either a polynomial of even degree or an exponential function.

We interpret the GFF $\varphi$ as the terminal value of a Brownian martingale $(W_t)_{t\in[0,1]}$ with values in $\mathcal{S}'(\mathbb{R}^2)$ and covariance given by

$$\mathbb{E}[W_t(x)W_s(y)] = (t \wedge s)(m^2 - \Delta)^{-1}(x-y), \qquad x,y \in \mathbb{R}^2, \quad t,s \in [0,1].$$

If we denote by $\mathbb{P}$ the law on the path space $C([0,1], \mathcal{S}'(\mathbb{R}^2))$ of this Brownian motion. The measure $\theta$ can be then interpreted as the law of $W_1$ under the Gibbs measure

$$\mathbb{Q}^\xi(dW) = \frac{1}{Z} e^{-V^\xi(W_1)} \mathbb{P}(dW).$$

It is well known that the measure $\mathbb{Q}^\xi$ is the minimizer of a Gibbs variational problem and that

$$-\log Z^\xi = -\log \mathbb{E}[e^{-V(W_1)}] = \inf_{\mathbb{R}} \{\mathbb{E}_\mathbb{R}[V(W_1)] + H(\mathbb{R}|\mathbb{P})\},$$

where $\mathbb{E}_\mathbb{R}$ denotes the expectation with respect to a probability measure $\mathbb{R}$ on the path space $C([0,1], \mathcal{S}'(\mathbb{R}^2))$ and $H(\mathbb{R}|\mathbb{P})$ the relative entropy of $\mathbb{R}$ with respect to $\mathbb{P}$. A formula of Boué–Dupuis [9, 35] allows to rewrite the variational problem as an optimal control problem, namely

$$-\log Z^\xi = -\log \mathbb{E}[e^{-V^\xi(W_1)}] = \inf_Z \mathbb{E}[V^\xi(W_1 + Z_1) + \mathcal{E}(Z)], \tag{1}$$

where the infimum is now taken over all processes $(Z_t)_{t\in[0,1]}$ with values in $H^1(\mathbb{R}^2)$ and adapted to the filtration of $W$ and

$$\mathcal{E}(Z) = \mathcal{E}(Z, Z) = \frac{1}{2} \int_0^1 \int_{\mathbb{R}^2} \dot{Z}_s (m^2 - \Delta) \dot{Z}_s ds$$

is a bilinear form giving the "energy" of the *drift* $Z$. Here $\dot{Z}_s$ denotes the derivative of $Z_s$ with respect to the parameter $s \in [0,1]$. As we show in this paper, the meaning of an optimal drift $Z$ is that of providing the law of the "interacting field" via the equality

$$\text{Law}_{\mathbb{Q}^\xi}(W) = \text{Law}_\mathbb{P}(W + Z) \tag{2}$$

valid for any minimizer $Z$ of the problem (1) and in particular

$$\theta = \mathrm{Law}_{\mathbb{P}}(W_1 + Z_1).$$

The existence of such minimizers will be guaranteed, in the rigorous discussion below, by suitably relaxing the variational problem into a problem on the space of probability measures. For pedagogical reasons let us for the moment continue our heuristic discussion ignoring these fine points.

The variational description (1) looses meaning in infinite volume, i.e. when $\xi \equiv 1$, since the cost functional is always infinite. However the main idea in this paper is to note that any minimizer $Z$ of (1) satisfy the Euler–Lagrange equation (first variation)

$$\mathbb{E}\big[\nabla V(W_1 + Z_1) \cdot K_1 + 2\mathcal{E}(Z, K)\big] = 0 \qquad (3)$$

for all "test controls" $K$, where $\nabla V(W_1 + Z_1) \cdot K_1$ denotes the derivative of $\varphi \mapsto V(\varphi)$ at $W_1 + Z_1$ and in the direction of $K_1$. In the following we will refer to the equation (3) as a *weak forward–backwards stochastic (partial) differential equation* (wFBSDE). The reason for this name derives from observing that (3) can be formally rewritten as

$$\int_0^1 \int_{\mathbb{R}^2} \mathbb{E}\big[(\nabla V(W_1 + Z_1) + (m^2 - \Delta)\dot{Z}_s)\dot{K}_s\big]\mathrm{d}s = 0$$

and since this holds for any nice adapted control $K$ it implies that

$$(m^2 - \Delta)\dot{Z}_s = -\mathbb{E}\big[\nabla V(W_1 + Z_1) \,|\, \mathcal{F}_s\big] \qquad (4)$$

which is a kind of elliptic equation for the optimal control $(Z_s)_{s \in [0,1]}$. It can be cast into a system of a forward SDE for $Z$ and a backward SDE which computes the conditional expectation. In this paper we will not exploit further this specific interpretation of (3) as a FBSDE so it is not necessary to make it more precise. We would like to stress however that our wFBSDE is a new kind of stochastic quantisation for Euclidean QFTs. Indeed, equation (3) makes sense also in infinite volume and can be used, very much like the weak formulation of an elliptic problem, to get estimates on the minimizer $Z$, e.g. by testing with $K = \rho Z$ where $\rho \colon \mathbb{R}^2 \to \mathbb{R}_+$ is some weight function decaying at spatial infinity. In particular we will obtain estimates of the form

$$\mathbb{E}\|\rho Z\|_{H_t^1 H_x^1}^2 < \infty,$$

uniformly in the cutoff $\xi$. This estimate together with the identification (2) is at the core of the proof of tightness of the family of measures $(\mathbb{Q}^\xi)_\xi$ with $\xi \to 1$. Any accumulation point $\mathbb{Q}$ represents an infinite volume Euclidean field, moreover the corresponding drift $Z$ gives, via the formula (2) an explicit coupling of $\mathbb{Q}$ and the GFF given by $\mathbb{P}$. The point of view which we advocate in this paper is that the wFBSDE equation (3) should provide a complete control over the drift $Z$ describing the infinite volume Euclidean field, for which a naive variational description is meaningless. Moreover intuitively we should be able to introduce a variational description of the measure $\mathbb{Q}$ with drift $Z$ via a renormalized BD formula of the form

$$-\log \mathbb{E}[e^{-f(X_1)}] = \inf_K \mathbb{E}\big[\,f(X_1 + K_1) + V(X_1 + K_1) - V(X_1) + \mathcal{E}(Z + K) - \mathcal{E}(Z)\,\big] \qquad (5)$$

with $X_s = W_s + Z_s$ being the "infinite volume interacting field" with potential $V$ and $f(\varphi)$ a test function which depends only locally on the field $\varphi$ and where the infimum is taken among all drifts $K$ which satisfy suitable decay at infinity. Given this decay one expects that the quantities $V(X_1 + K_1) - V(X_1)$ and $\mathcal{E}(Z + K) - \mathcal{E}(Z)$ are finite even if each summand is infinite.

Leaving this heuristic account of the ideas investigated in this paper, we describe now the plan for the rigorous proofs. We will start in Section 1 by introducing the general setting for the variational description of two dimensional Euclidean quantum fields as outlined above. We prove certain useful bounds and also that the minimizers of this variational problem describe in a precise sense the Euclidean field and moreover that they satisfy wFBSDEs.

As a first example, we address in Section 2 the problem of the infinite volume limit of the $\varphi_2^4$ measure $\theta_\Lambda$ on a bounded domain $\Lambda \subseteq \mathbb{R}^2$. From the wFBSDE we will derive a-priori bounds for minimizers strong enough to give tightness of the family $(\theta_\Lambda)_\Lambda$ as $\Lambda \nearrow \mathbb{R}^2$. Unfortunately, establishing rigorously the form (5) of the limiting variational problem requires certain locality properties of the solutions to the wFBSDE (3) in infinite volume. This is usually possible when $V$ is convex. Ultraviolet renormalization in $\varphi_2^4$ however spoils any "easy" convexity of the potential and it is not clear how to gain sufficient control of (3) in this case. Similar issues appear in other methods of stochastic quantization as soon as the infinite volume limit is considered for $P(\varphi)_2$ Euclidean quantum fields [18, 19, 1]. Despite these difficulties, in Section 3 we show that the wFBSDE is still a useful tool which can be used effectively to prove the large deviations of the $\varphi_2^4$ measure in the semiclassical limits. Let us remark that Lacoin Rhodes and Vargas obtained large deviations results in the case of the Liouville model in [26, 27].

Finally, in Section 4, we discuss the case of the exponential interaction, that is the $\exp(\beta\varphi)_2$ model (or Høegh–Krohn model) introduced in [21]. See e.g. [1] for a discussion of its stochastic quantization in the elliptic setting. For this model we repeat the construction of the variational description of the infinite volume measure. At variance with what happens in $P(\varphi)_2$ theories, renormalization does not spoil the convexity properties of the potential and therefore the study of the related wFBSDE can go further and in particular we will be able to prove the necessary locality of the perturbed equation. This in turn will allow us to prove the $\Gamma$-convergence of the variational problem along with the infinite volume limit and establish eq. (5) in its rigorous form. As a result of this analysis, we are able to obtain a concise and well-behaved description of an infinite volume Gibbs measure related to an Euclidean quantum field.

**Acknowledgments.** MG would like to thank Ronan Herry for pointing out the reference [29] on the Brownian transport map. MG has been partially supported by the Deutsche Forschungsgemeinschaft (DFG, German Research Foundation) through the Hausdorff Center for Mathematics under Germany's Excellence Strategy – EXC-2047/1 – 390685813 and through CRC 1060 - project number 211504053. NB is supported by the ERC Advanced Grant 741487 "Quantum Fields and Probability".

**Notation.** Let $\langle \alpha \rangle := (1+|\alpha|^2)^{1/2}$. Let $\mathcal{S}(\mathbb{R}^2)$ the Schwarz space of test functions on $\mathbb{R}^2$ and $\mathcal{S}'(\mathbb{R}^2)$ the corresponding space of tempered distributions. For any $\alpha \in \mathbb{R}$ we define fractional derivatives $\langle D \rangle^\alpha = (1+D^2)^{\alpha/2}$ on $\mathcal{S}'(\mathbb{R}^2)$ by

$$\langle D \rangle^\alpha f = \mathcal{F}^{-1}(\xi \mapsto \langle \xi \rangle^{\alpha/2} (\mathcal{F}f)(\xi)), \qquad f \in \mathcal{S}'(\mathbb{R}^2)$$

where $\mathcal{F}$ is the Fourier transform on $\mathbb{R}^2$.

All along the paper we need to consider Banach subspaces of $\mathcal{S}'(\mathbb{R}^2)$ associated with a positive spatial weight $\rho$. Let $\rho \colon \mathbb{R}^2 \to \mathbb{R}_{>0}$ a smooth positive function such that for some sufficiently small $\varepsilon > 0$ we have

$$|\Delta \rho^{1/2}| + |\nabla \rho^{1/2}| \leqslant \varepsilon \rho^{1/2}, \qquad \rho \lesssim \langle \cdot \rangle^{-32}. \tag{6}$$

The number $\varepsilon$ can depend on various parameters but not on the volume cutoffs we introduce to regularize the Euclidean QFTs. For some $\alpha > 0$, $s \in \mathbb{R}$, $p, q \in [1, \infty]$, we introduce the usual weighted Sobolev spaces $W^{s,p}(\rho^\alpha)$ and weighted Besov spaces $B^s_{p,q}(\rho^\alpha)$ with norms

$$\|f\|_{W^{s,p}(\rho^\alpha)} = \|\rho^\alpha \langle D \rangle^s f\|_{L^p(\mathbb{R}^2)}, \qquad \|f\|_{B^s_{p,q}(\rho^\alpha)} = \|\rho^\alpha f\|_{B^s_{p,q}},$$

and let $H^s(\rho^\alpha) = W^{1,2}(\rho^\alpha)$ and $C^s(\rho^\alpha) = B^s_{\infty,\infty}(\rho^\alpha)$. Here $B^s_{p,q}$ is the unweighted Besov norm given by

$$\|f\|_{B^s_{p,q}} = \sum_{i \geqslant -1} 2^{isq} \|\Delta_i f\|_{L^p(\mathbb{R}^2)},$$

where the Littlewood–Paley operators $(\Delta_i)_{i \geqslant -1}$ are defined by $\Delta_i = \psi(D/2^i)$ for $i \geqslant 0$ and some function $\psi : \mathbb{R}^2 \to \mathbb{R}_{\geqslant 0}$ supported in an annulus centered at the origin and $\Delta_{-1} = \tilde{\psi}(D)$ for $\tilde{\psi} : \mathbb{R}^2 \to \mathbb{R}_{\geqslant 0}$ supported on a ball centered at the origin and such that $\sum_{i \geqslant -1} \Delta_i = 1$.

In general we avoid to indicate the weight when dealing with the unweighted ($\alpha = 0$) version of the above spaces.

We denote by $\mathcal{P}(X)$ the set of Borel measures on the topological space $X$, for $P \in \mathcal{P}(X)$ we denote $\mathbb{E}_P$ the corresponding expectation and for any random variable $H : X \to Y$ on the topological space $Y$ we let $\text{Law}_P(H) \in \mathcal{P}(Y)$ the corresponding law.

# 1  A variational setting for 2d Euclidean fields

Consider the Gibbs measure

$$\theta^V(\mathrm{d}\phi) := \frac{e^{-V(\phi)} \theta^0(\mathrm{d}\phi)}{\int e^{-V(\phi)} \theta^0(\mathrm{d}\phi)}, \qquad \phi \in \mathcal{S}'(\mathbb{R}^2), \tag{7}$$

where $\theta^0$ is the law of the Gaussian field $\phi$ (GFF) on $\mathbb{R}^2$ with covariance $(m^2 - \Delta)^{-1}$ for some mass $m > 0$ (fixed throughout the paper) and where $\Delta$ denotes the Laplacian on $\mathbb{R}^2$. The function $V : \mathcal{S}'(\mathbb{R}^2) \to \mathbb{R}$ represent an interaction potential and for now we will only assume that it is measurable and *tame* with respect to $\theta^0$ according to the following definition.

**Definition 1.** *Let $(X, \mathcal{F}, \gamma)$ be a probability space. A measurable function $F : X \to [-\infty, \infty]$ is called tame with respect to $\gamma$ if there exists $p, q > 1$ with $\frac{1}{p} + \frac{1}{q} = 1$ and*

$$\mathbb{E}_\gamma[\exp(-qF)] + \mathbb{E}_\gamma[|F|^p] < \infty.$$

Tameness has been introduced by Üstunel [35] as a convenient condition for the validity of the Boué–Dupuis variational formula for Brownian functionals. It has been recently improved by Hariya and Watanabe [20] to allow for $p = q = 1$. To introduce this formula we need to enlarge our probability space and realize the GFF as a marginal of a Brownian motion, following the approach initiated in [3]. Here we are able to use a simpler embedding thanks to the milder nature of the Gaussian free field in two dimensions compared to its three dimensional counterpart.

Let $\mathbb{P}$ be the law of a cylindrical Brownian motion $(X_t)_{t \in [0,1]}$ on $L^2(\mathbb{R}^2)$. It can be realized on $C([0,1], C^{-1-\varepsilon}(\rho^\varepsilon))$ where $\varepsilon > 0$ can be taken arbitrarily small. Let

$$W_t := (m^2 - \Delta)^{-1/2} X_t, \qquad t \in [0,1],$$

and note that indeed $\theta^0 = \mathrm{Law}_\mathbb{P}(W_1)$. We will consider $W$ as a random variable in $\mathfrak{S}_0 := C([0,1], C^{-\varepsilon}(\rho^\varepsilon))$.

Let $\mathfrak{H}(\rho^\alpha) = H_t^1 H_x^1(\rho^\alpha)$, $\mathfrak{H} = \mathfrak{H}(\rho^0)$. Denote with $\mathfrak{H}^a(\rho^\alpha)$ the space of all processes $u: \mathfrak{S}_0 \times [0,1] \to \mathbb{R}$ progressively measurable with respect to the filtration generated by $(W_t)_{t \in [0,1]}$, starting from 0 and equipped with the norm of $L^2(\mathbb{P}, \mathfrak{H}(\rho^\alpha))$. Moreover define $\mathcal{E}: \mathfrak{H}^a \to L^0(\mathbb{P})$ as

$$\mathcal{E}(Z) := \frac{1}{2} \int_0^1 \|(m^2 - \Delta)^{1/2} \dot{Z}_s\|_{L^2}^2 \mathrm{d}s = \frac{1}{2} \int_0^1 (\|\nabla \dot{Z}_s\|_{L^2}^2 + m^2 \|\dot{Z}_s\|_{L^2}^2) \mathrm{d}s$$

where $\dot{Z}$ denotes the time-derivative of $Z$. With these definitions, the Boué–Dupuis formula reads as follows.

**Lemma 2.** *(Boué–Dupuis formula)* For any tame functional $G: \mathfrak{S}_0 \to \mathbb{R}$ we have

$$-\log \mathbb{E}[e^{-G(W_\cdot)}] = \inf_{Z \in \mathfrak{H}^a} \mathbb{E}[G(W_\cdot + Z_\cdot) + \mathcal{E}(Z)].$$

**Proof.** Theorem 7 from Üstunel [35] specialized in our context (see also [3]) gives

$$-\log \mathbb{E}[e^{-G(W_\cdot)}] = \inf_{u \in \mathbb{H}^a} \mathbb{E}\left[G(W_\cdot + Z_\cdot) + \frac{1}{2} \int_0^1 \|u_s\|_{L^2}^2 \mathrm{d}s\right],$$

where $\mathbb{H}^a$ is the space of all progressively measurable processes starting from 0 equipped with the norm of $L^2(\mathbb{P}, L_t^2 L_x^2)$ and

$$Z_t := (m^2 - \Delta)^{-1/2} \int_0^t u_s \mathrm{d}s.$$

Observing that $u_t = (m^2 - \Delta)^{1/2} \dot{Z}_t$ we obtain the claim. $\square$

Note that the functional $G(W_\cdot + Z_\cdot)$ makes sense for any $Z \in \mathfrak{H}^a$ since by Theorem 2 from Üstunel [35] we have that the push-forward measure $\mathbb{Q} = (W_\cdot + Z_\cdot)_* \mathbb{P}$ has finite relative entropy with respect to $\mathbb{P}$ and therefore is absolutely continuous, i.e. $\mathbb{Q} \ll \mathbb{P}$.

We call a random variable with values in $\mathfrak{H}^a$ a *control* (or *drift*) and we will understand that $W = W_1$ and $Z = Z_1$ whenever this does not lead to confusion.

To simplify certain technical considerations later is convenient to augment the process $W$ with some additional informations. Consider a Polish space $\mathfrak{S}$ endowed with a continuous projection $\pi_W: \mathfrak{S} \to \mathfrak{S}_0$ and let $\mathbb{W}: \mathfrak{S}_0 \to \mathfrak{S}$ a measurable map such that $\pi_W \mathbb{W} = W$, i.e. $\mathfrak{S} \xrightarrow{W} \mathfrak{S}_0$ is a fibered topological space and $\mathbb{W}$ a measurable section. We will understand $\mathbb{W}$ as an *enhanced GFF*. The kind of enhancements needed will be different in different models. We will always consider the space $\mathfrak{S}$ endowed with the measure given by the law of $\mathbb{W}$ under $\mathbb{P}$. A functional on $\mathfrak{S}$ is tame if it is tame for this measure.

For any measurable $f: \mathcal{S}'(\mathbb{R}^2) \to \mathbb{R}$ bounded below and in $L^2(\theta^0)$ (not necessarily linear) define $\mathcal{W}^V(f)$ as

$$e^{-\mathcal{W}^V(f)} := \int e^{-f(\phi)} \theta^V(\mathrm{d}\phi).$$

As a direct consequence of Lemma 2, the functional $\mathcal{W}^V$ has the variational representation

$$\mathcal{W}^V(f) = \inf_{Z \in \mathfrak{H}^a} F^{f+V}(Z) - \inf_{Z \in \mathfrak{H}^a} F^V(Z).$$

for any $f$ such that $f+V$ is tame and where for a tame functional $G$ on $\mathfrak{S}_0$ we let

$$F^G(Z) := \mathbb{E}[G(W+Z) + \mathcal{E}(Z)], \qquad Z \in \mathfrak{H}^a.$$

To proceed we want to find a minimizer of $F^{f+V}$. It is however not clear if we can find one in $\mathfrak{H}^a$. Following the approach we used in [3], we relax our problem and instead minimize over the law of the pair $(W, Z)$ with $W$ fixed and $Z$ varying over $\mathfrak{H}^a$. In order to ensure compactness we need to complete the space of trial laws for a suitable topology for which our target functional is coercive.

Let $\mathfrak{H}_w(\rho^\alpha)$ to be $\mathfrak{H}(\rho^\alpha)$ equipped with the weak topology and let $\hat{\mathfrak{H}} := \mathfrak{H}(\rho^{1/2})$, $\hat{\mathfrak{H}}_w := \mathfrak{H}_w(\rho^{1/2})$ and $\hat{\mathfrak{H}}^a := \mathfrak{H}^a(\rho^{1/2})$. With a slight abuse of notation we denote $(\mathbb{W}, Z)$ the canonical coordinates on $\mathfrak{S} \times \hat{\mathfrak{H}}_w$.

**Definition 3.** *Let*

$$\mathcal{X} := \left\{ \mu = \mathrm{Law}_{\mathbb{P}}(\mathbb{W}, Z) \in \mathcal{P}(\mathfrak{S} \times \hat{\mathfrak{H}}_w) : Z \in \hat{\mathfrak{H}}^a, \mathbb{E}_\mu[\|Z\|_{\hat{\mathfrak{H}}}^2] < \infty \right\},$$

*and*

$$\bar{\mathcal{X}} := \left\{ \mu \in \mathcal{P}(\mathfrak{S} \times \hat{\mathfrak{H}}_w) : \exists \mu_n \in \mathcal{X}, \mu_n \to \mu \text{ weakly}, \sup_{n \in \mathbb{N}} \mathbb{E}_{\mu_n}[\|Z\|_{\hat{\mathfrak{H}}}^2] < \infty \right\}.$$

*We equip $\bar{\mathcal{X}}$ with the following notion of convergence: a sequence $(\mu_n)_n$ in $\bar{\mathcal{X}}$ converges to $\mu$ if*

a) $\mu_n \to \mu$ *weakly on* $\mathfrak{S} \times \hat{\mathfrak{H}}_w$,

b) $\sup_{n \in \mathbb{N}} \mathbb{E}_{\mu_n}[\|Z\|_{\hat{\mathfrak{H}}}^2] < \infty$.

*Note that the probability space in the definition of $\mathcal{X}$ is fixed, so all elements in $\mathcal{X}$ are laws of random variables living on a fixed probability space. Note that all elements of $\mathcal{X}$ have first marginal given by $\mathrm{Law}(\mathbb{W})$, and this also carries over to $\bar{\mathcal{X}}$ by weak convergence.*

We will need to realize couplings of two elements of $\bar{\mathcal{X}}$. These are constructed by the next lemma which needs the following preliminary definition of the space hosting the coupling. Denote the canonical coordinates on $\mathfrak{S} \times \hat{\mathfrak{H}}_w \times \hat{\mathfrak{H}}_w$ by $(\mathbb{W}, Z, K)$.

**Definition 4.** *For $\mu \in \bar{\mathcal{X}}$ let $\bar{\mathcal{Y}}(\mu) \subseteq \mathcal{P}(\mathfrak{S} \times \hat{\mathfrak{H}}_w \times \hat{\mathfrak{H}}_w)$ the set of all $\upsilon$ such that $\mathrm{Law}_\upsilon(\mathbb{W}, Z) = \mu$ and for which there exists sequences $(Z^n)_n, (K^n)_n \subseteq \hat{\mathfrak{H}}^a$ with*

$$\mathrm{Law}_{\mathbb{P}}(\mathbb{W}, Z^n, K^n) \to \upsilon \text{ and } \sup_{n \in \mathbb{N}} \mathbb{E}[\|Z^n\|_{\hat{\mathfrak{H}}}^2 + \|K^n\|_{\hat{\mathfrak{H}}}^2] < \infty.$$

**Lemma 5.** *For any $\nu, \mu \in \bar{\mathcal{X}}$ there exists $\upsilon \in \bar{\mathcal{Y}}(\mu)$ such that*

$$\mathrm{Law}_\upsilon(\mathbb{W}, Z+K) = \mathrm{Law}_\nu(\mathbb{W}, Z).$$

**Proof.** Since $\mu, \nu$ are in $\bar{\mathcal{X}}$ we can find sequences $Z^n, \bar{Z}^n \in \hat{\mathfrak{H}}^a$ such that $\mu_n := \mathrm{Law}_{\mathbb{P}}(\mathbb{W}, Z^n) \to \mu$ and $\nu_n := \mathrm{Law}_{\mathbb{P}}(\mathbb{W}, \bar{Z}^n) \to \nu$ weakly in $\hat{\mathfrak{H}}_w$ and furthermore $\sup_n \mathbb{E}[\|Z^n\|_{\hat{\mathfrak{H}}}^2 + \|\bar{Z}^n\|_{\hat{\mathfrak{H}}}^2] < \infty$. Now define

$$K^n = \bar{Z}^n - Z^n,$$

and set

$$\upsilon^n = \mathrm{Law}_{\mathbb{P}}(\mathbb{W}, Z^n, K^n).$$

We have that $\sup_n \mathbb{E}[\|K^n\|_{\hat{\mathfrak{H}}}^2] \leqslant \sup_n \mathbb{E}[\|Z^n\|_{\hat{\mathfrak{H}}}^2 + \|\bar{Z}^n\|_{\hat{\mathfrak{H}}}^2] < \infty$ so $\upsilon^n$ is tight on $\mathfrak{S} \times \hat{\mathfrak{H}}_w \times \hat{\mathfrak{H}}_w$. Taking a convergent subsequence (not relabeled) we can assume that $\upsilon^n \to \upsilon$. We then have

$$\begin{aligned}
\text{Law}_\upsilon(\mathbb{W}, Z+K) &= \lim_{n\to\infty} \text{Law}_{\upsilon^n}(\mathbb{W}, Z+K) \\
&= \lim_{n\to\infty} \text{Law}_\mathbb{P}(\mathbb{W}, Z^n + K^n) \\
&= \lim_{n\to\infty} \text{Law}_\mathbb{P}(\mathbb{W}, \bar{Z}^n) = \nu
\end{aligned}$$

which proves the claim. $\square$

On functionals on $\bar{\mathcal{X}}$ we can define first and second order derivatives as follows. Let $\mu \in \bar{\mathcal{X}}$ and let $K: \mathfrak{S} \times \hat{\mathfrak{H}}_w \to \hat{\mathfrak{H}}_w$ a measurable process adapted to the filtration generated by $(\mathbb{W}_t, Z_t)_{t \in [0,1]}$, continuous in the second component and such that $K \in L^2(\mu, \hat{\mathfrak{H}})$. We call such $K$ an admissible control for $\mu$. For all $\sigma \in \mathbb{R}$ let $\mu^{\sigma K}$ the law of $(\mathbb{W}, Z + \sigma K)$ under $\mu$. We have $\mu^{\sigma K} \in \bar{\mathcal{X}}$. A functional $\check{F}: \bar{\mathcal{X}} \to \mathbb{R}$ is differentiable in $\bar{\mathcal{X}}$ in the direction of the admissible control $K$ if the map $\sigma \mapsto \check{F}(\mu^{\sigma K})$ is differentiable in $\sigma = 0$ and in this case we define

$$\nabla \check{F}(\mu)(K) = \frac{\mathrm{d}}{\mathrm{d}\sigma}\bigg|_{\sigma=0} \check{F}(\mu^{\sigma K}).$$

Second order derivatives can be similarly defined by taking two admissible controls $K, K'$ and letting

$$\nabla^2 \check{F}(\mu)(K, K') = \frac{\mathrm{d}}{\mathrm{d}\sigma'}\bigg|_{\sigma'=0} \frac{\mathrm{d}}{\mathrm{d}\sigma}\bigg|_{\sigma=0} \check{F}(\mu^{\sigma K + \sigma' K'}),$$

whenever $(\sigma, \sigma') \mapsto \check{F}(\mu^{\sigma K + \sigma' K})$ is $C^2$ in $(\sigma, \sigma') = 0$.

We are ready to relax the variational problem and prove an alternative representation of $\mathcal{W}_\Lambda(f)$. For a tame $G: \mathfrak{S} \to \mathbb{R}$ define the functional $\check{F}^G: \bar{\mathcal{X}} \to \mathbb{R}$ by

$$\check{F}^G(\mu) := \mathbb{E}_\mu[G(W+Z) + \mathcal{E}(Z)], \qquad \mu \in \bar{\mathcal{X}}.$$

**Definition 6.** *We say that a functional $\check{F}: \bar{\mathcal{X}} \to \mathbb{R}$ is admissible if it is lower semicontinuous and such that*

$$\check{F}(\mu) \gtrsim 1 + \mathbb{E}_\mu[\mathcal{E}(Z)], \tag{8}$$

*this ensures in particular that the functional $\check{F}: \bar{\mathcal{X}} \to \mathbb{R}$ is coercive, that is for all $L < \infty$ the level sets $\{\check{F} \leqslant L\}$ are compact in $\bar{\mathcal{X}}$.*

**Definition 7.** *For a tame $G: \mathfrak{S} \to \mathbb{R}$ define the functional $\check{F}^G: \bar{\mathcal{X}} \to \mathbb{R}$ by*

$$\check{F}^G(\mu) := \mathbb{E}_\mu[G(W+Z) + \mathcal{E}(Z)], \qquad \mu \in \bar{\mathcal{X}}.$$

*We say that $G$ is admissible if $G$ is tame, $\check{F}^G$ is admissible and $(W, Z) \to G(W+Z)$ is continuous on $\mathfrak{S} \times \hat{\mathfrak{H}}_w$.*

The following lemma is an obvious consequence of this relaxed variational setting.

**Lemma 8.** *Assume that $G: \mathfrak{S} \to \mathbb{R}$ is admissible, then*

$$\inf_{Z \in \mathfrak{H}^a} F^G(Z) = \inf_{Z \in \hat{\mathfrak{H}}^a} F^G(Z) = \inf_{\mu \in \mathcal{X}} \check{F}^G(\mu) = \min_{\mu \in \bar{\mathcal{X}}} \check{F}^G(\mu). \tag{9}$$

Let $\mu^G \in \bar{\mathcal{X}}$ be any minimizer of $\check{F}^G$. If $\check{F}^G$ is differentiable at $\mu^G$ it satisfies the equation $\nabla \check{F}^G(\mu^G)(K) = 0$ for all admissible controls $K$. Moreover if it is twice differentiable at $\mu^G$ we also have the positivity condition $\nabla^2 \check{F}^G(\mu^G)(K, K') \geqslant 0$ for all admissible controls $K$, $K'$.

**Proof.** The equality $\inf_{Z \in \mathfrak{H}^a} F^G(Z) = \inf_{Z \in \hat{\mathfrak{H}}^a} F^G(Z)$ follows from the property (8) of the admissible functional $\check{F}^G$ since $\mathfrak{H}^a$ can be viewed as a subspace of $\bar{\mathcal{X}}$ via the map $u \mapsto \mathrm{Law}(\mathbb{W}, u)$. The equality $\inf_{Z \in \hat{\mathfrak{H}}^a} F^G(Z) = \inf_{\mu \in \mathcal{X}} \check{F}^G(\mu)$ follows directly from the definition of $\mathcal{X}$. We have that $\inf_{\mu \in \bar{\mathcal{X}}} \check{F}^G(\mu) = \min_{\mu \in \bar{\mathcal{X}}} \check{F}^G(\mu)$ again from the admissibility of $G$ which ensures that sequences $(\mu_n)_n \subseteq \mathcal{X}$ of approximate minimizers of $\check{F}^G$ converge in $\bar{\mathcal{X}}$ and that any limit $\mu^G$ is a minimizer. The first and second order conditions are then consequences of the differentiability of $\check{F}^G$. It remains to show

$$\inf_{\mu \in \mathcal{X}} \check{F}^G(\mu) = \inf_{\mu \in \bar{\mathcal{X}}} \check{F}^G(\mu).$$

We fist show this statement for bounded $G$. We want to show that for any $\mu$ in $\bar{\mathcal{X}}$ any $\delta > 0$ there exists $\mu_n \in \mathcal{X}$

$$\mathbb{E}_{\mu_n}[G(W+Z) + \mathcal{E}(Z)] \leqslant \mathbb{E}_\mu[G(W+Z) + \mathcal{E}(Z)] + \delta.$$

So given $\mu \in \bar{\mathcal{X}}$ we take a sequence $\mu_n \in \mathcal{X}$ such that $\mu_n \to \mu$. We define $\mu_n^N$ to be

$$\mathrm{Law}(\mathbb{W}, Z^N) \quad \text{where } Z_t^N = Z_{t \wedge \tau_N} \quad \text{where } \tau_N = \inf_{t \in [0,1]} \int_0^t \|\nabla \dot{Z}_s\|_{L^2}^2 \mathrm{d}s \geqslant N.$$

Then $\mu_n^N = \mu_n$ on the event that $\mathcal{E}(Z) \leqslant N$. Select a subsequence such that $\mu_n^N \to \upsilon^N$ in $\mathcal{X}$. Define $\theta_*^\delta \nu = \mathrm{Law}_\nu(W, \theta^\delta * Z)$ where $\theta^\delta$ is a sequence of standard mollifiers. Now by weak convergence

$$\lim_{n \to \infty} \mathbb{E}_{\theta_*^\delta \mu_n^N}[G(W+Z) + \mathcal{E}(Z)] = \mathbb{E}_{\theta_*^\delta \upsilon^N}[G(W+Z) + \mathcal{E}(Z)],$$

since the functional inside the expectation is continuous and bounded on the joint support of the measures $\theta_*^\delta \mu_n^N$. We have to show that

$$\mathbb{E}_{\theta_*^\delta \upsilon^N}[G(W+Z) + \mathcal{E}(Z)] \leqslant \mathbb{E}_\mu[G(W+Z) + \mathcal{E}(Z)] + \bar{\delta}.$$

First we establish that

$$\lim_{\delta \to 0} \mathbb{E}_{\theta_*^\delta \upsilon^N}[G(W+Z) + \mathcal{E}(Z)] = \mathbb{E}_{\upsilon^N}[G(W+Z) + \mathcal{E}(Z)]$$

Since $\mathcal{E}(\theta^\delta * Z) \leqslant \mathcal{E}(Z)$ and $G$ is continuous on $\mathfrak{S} \times \hat{\mathfrak{H}}_w$ we have

$$\lim_{\delta \to 0} \mathbb{E}_{\theta_*^\delta \upsilon^N}[G(W+Z)] = \mathbb{E}_{\upsilon^N}[G(W+Z)] \qquad \lim_{\delta \to 0} \mathbb{E}_{\upsilon^N}[\mathcal{E}(\theta^\delta * Z)] = \mathbb{E}[\mathcal{E}(Z)].$$

Now we have to show that

$$\mathbb{E}_{\upsilon^N}[G(W+Z) + \mathcal{E}(Z)] \leqslant \mathbb{E}_\mu[G(W+Z) + \mathcal{E}(Z)] + \delta.$$

Since $\upsilon^N$ coincides with $\mu$ on any open set contained in $\mathcal{E}(Z) \leqslant N$ we have

$$|\mathbb{E}_{\upsilon^N}[G(W+Z)] - \mathbb{E}_\mu[G(W+Z)]| \leqslant \|G\|_{L^\infty} \mathbb{E}_\mu\big[\mathbb{1}_{\mathcal{E}(Z) > N/2}\big] \lesssim N^{-1}.$$

For the energy

$$\begin{aligned} \mathbb{E}_{\upsilon^N}[\mathcal{E}(Z)] &\leqslant \liminf_{n \to \infty} \mathbb{E}_{\mu_n^N}[\mathcal{E}(Z)] \leqslant \liminf_{n \to \infty} \mathbb{E}_{\mu_n}\left[\int_0^{\tau_N} \|\nabla \dot{Z}_s^\delta\|_{L^2}^2 \mathrm{d}s\right] \\ &\leqslant \liminf_{n \to \infty} \mathbb{E}_{\mu_n}[\mathcal{E}(Z) \wedge N] \leqslant \mathbb{E}_\mu[\mathcal{E}(Z) \wedge N]. \end{aligned}$$

Now
$$\mathbb{E}_\mu[\mathcal{E}(Z) \wedge N] \leqslant \mathbb{E}_\mu[\mathcal{E}(Z)].$$

We now remove the restriction that $G$ is bounded. To this end let $G^N = (G \wedge N) \vee (-N)$. Then we have established that already that
$$-\log \int \exp(-G^N) \mathrm{d}\theta^0 = \inf_{\mu \in \mathcal{X}} \check{F}^{G^N}(\mu) = \min_{\mu \in \bar{\mathcal{X}}} \check{F}^{G^N}(\mu).$$

Now let $\nu$ be a minimizer of $\check{F}^G$ in $\bar{\mathcal{X}}$ (existence of this was established above). Then

$$\begin{aligned}
\inf_{\mu \in \bar{\mathcal{X}}} \mathbb{E}_\mu[G(W+Z) + \mathcal{E}(Z)] &= \mathbb{E}_\nu[G(W+Z) + \mathcal{E}(Z)] \\
&= \lim_{N \to \infty} \mathbb{E}_\nu[G^N(W+Z) + \mathcal{E}(Z)] \\
&\geqslant \lim_{N \to \infty} \inf_{\mu \in \bar{\mathcal{X}}} \mathbb{E}_\mu[G^N(W+Z) + \mathcal{E}(Z)] \\
&= \lim_{N \to \infty} -\log \int \exp(-G^N) \mathrm{d}\theta_0 \\
&= -\log \int \exp(-G) \mathrm{d}\theta_0 \\
&= \inf_{\mu \in \mathcal{X}} \mathbb{E}_\mu[G(W+Z) + \mathcal{E}(Z)]
\end{aligned}$$

Where we have used dominated convergence. $\square$

**Remark 9.** Note that
$$\check{F}^G(\mu^{\sigma K}) - \check{F}^G(\mu) = \mathbb{E}_\mu[G(W+Z+\sigma K) - G(W+Z) + \mathcal{E}(Z+\sigma K) - \mathcal{E}(Z)]$$
$$= \mathbb{E}_\mu[G(W+Z+\sigma K) - G(W+Z) + 2\sigma \mathcal{E}(Z,K) + \sigma^2 \mathcal{E}(K)].$$

If $G$ is differentiable in the direction of $H^1$ we have
$$\nabla \check{F}^G(\mu)(K) = \mathbb{E}_\mu[\nabla G(W+Z) \cdot K + 2\mathcal{E}(Z,K)].$$

If $G$ is twice differentiable we have also
$$\nabla^2 \check{F}^G(\mu)(K,K') = \mathbb{E}_\mu[\nabla^2 G(W+Z) \cdot (K,K') + 2\mathcal{E}(K,K')].$$

The consequences of the existence of this minimizer are clarified by the following lemmas.

**Lemma 10.** *If $f+V$ is admissible and $\mu^{f+V}$ a minimizer of $\check{F}^{f+V}$ then we have*
$$\mathbb{E}_{\mu^{f+V}}[f(W+Z)] \leqslant \mathcal{W}^V(f) \leqslant \mathbb{E}_{\mu^V}[f(W+Z)]. \tag{10}$$

**Proof.** If $f+V$ is admissible we have proven that
$$\mathcal{W}^V(f) = \inf_{\mu \in \bar{\mathcal{X}}} \check{F}^{f+V}(\mu) - \inf_{\mu \in \bar{\mathcal{X}}} \check{F}^V(\mu) = \check{F}^{f+V}(\mu^{f+V}) - \check{F}^V(\mu^V). \tag{11}$$

Eq. (11) implies in turn
$$\check{F}^{f+V}(\mu^{f+V}) - \check{F}^V(\mu^{f+V}) \leqslant \mathcal{W}^V(f) \leqslant \check{F}^{f+V}(\mu^V) - \check{F}^V(\mu^V)$$

and since
$$\check{F}^{f+V}(\mu^{g+V}) - \check{F}^V(\mu^{g+V}) = \mathbb{E}_{\mu^{g+V}}[f(W+Z) + V(W+Z) + \mathcal{E}(Z)]$$
$$-\mathbb{E}_{\mu^{g+V}}[V(W+Z) + \mathcal{E}(Z)] = \mathbb{E}_{\mu^{g+V}}[f(W+Z)]$$

for any $g$ (in particular admissible $g$), we obtain eq. (10). $\square$

**Lemma 11.** *Let $V\colon \mathfrak{S} \to \mathbb{R}$ be admissible. Then for all $f\colon \mathfrak{S} \to \mathbb{R}$ continuous and bounded $f+V$ is admissible. Let $\mu^{f+V}$ be a minimizer for $\check{F}^{f+V}$, then for all $\alpha \in \mathbb{R}$ we have*

$$\frac{\mathrm{d}}{\mathrm{d}\alpha}\mathcal{W}^V(\alpha f) = \mathbb{E}_{\mu^{\alpha f + V}}[f(W+Z)]. \tag{12}$$

**Proof.** The admissibility of $f+V$ is easy to obtain given that $f$ is bounded and continuous and $V$ is admissible. It is clear that $\alpha \mapsto \mathcal{W}^V(\alpha f)$ is differentiable for bounded $f$, we need only to establish the representation of this derivative via the variational minimizer $\mu^{\alpha f}$. For any $\gamma \geqslant 0$ we have, by direct computation

$$\mathcal{W}^V((\alpha+\gamma)f) - \mathcal{W}^V(\alpha f) \leqslant \check{F}^{(\alpha+\gamma)f+V}(\mu^{\alpha f+V}) - \check{F}^{\alpha f+V}(\mu^{\alpha f+V}) = \gamma \mathbb{E}_{\mu^{\alpha f+V}}[f(W+Z)],$$

and

$$\mathcal{W}^V(\alpha f) - \mathcal{W}^V((\alpha-\gamma)f) \geqslant \check{F}^{\alpha f+V}(\mu^{\alpha f+V}) - \check{F}^{(\alpha-\gamma)f+V}(\mu^{\alpha f+V}) = \gamma \mathbb{E}_{\mu^{\alpha f+V}}[f(W+Z)].$$

By differentiability of $\mathcal{W}^V(\alpha f)$ we have

$$\mathbb{E}_{\mu^{\alpha f+V}}[f(W+Z)] \leqslant \lim_{\gamma \to 0+} \frac{\mathcal{W}^V((\alpha+\gamma)f) - \mathcal{W}^V(\alpha f)}{\gamma} = \frac{\mathrm{d}}{\mathrm{d}\alpha}\mathcal{W}^V(\alpha f)$$

and

$$\frac{\mathrm{d}}{\mathrm{d}\alpha}\mathcal{W}^V(\alpha f) = \lim_{\gamma \to 0+} \frac{\mathcal{W}^V(\alpha f) - \mathcal{W}^V((\alpha-\gamma)f)}{\gamma} \leqslant \mathbb{E}_{\mu^{\alpha f+V}}[f(W+Z)]$$

from which we can conclude. $\square$

There are several interesting consequences of this formula. First of all, integrating (12) in $\alpha$, we have

$$\mathcal{W}^V(f) = \int_0^1 \frac{\mathrm{d}}{\mathrm{d}\alpha}\mathcal{W}^V(\alpha f)\mathrm{d}\alpha = \int_0^1 \mathbb{E}_{\mu^{\alpha f+V}}[f(W+Z)]\mathrm{d}\alpha, \tag{13}$$

which is a representation of $\mathcal{W}^V(f)$ alternative to the variational formula. On the other hand, taking $\alpha = 0$ in (12) we obtain

$$\int f(\phi)\theta^V(\mathrm{d}\phi) = \frac{\mathrm{d}}{\mathrm{d}\alpha}\bigg|_{\alpha=0}\mathcal{W}^V(\alpha f) = \mathbb{E}_{\mu^V}[f(W_1+Z_1)].$$

**Theorem 12.** *The law of $W_1+Z_1$ under $\mu^V$ is indeed $\theta^V$. In particular for all $f\colon \mathfrak{S} \to \mathbb{R}$ continuous and bounded*

$$\mathcal{W}^V(f) = -\log \int e^{-f(\phi)}\theta^V(\mathrm{d}\phi) = -\log \mathbb{E}_{\mu^V}[e^{-f(W_1+Z_1)}].$$

This last formula naturally calls for a variational description of the average $\mathbb{E}_{\mu^V}[e^{-f(W_1+Z_1)}]$ with respect to the base measure $\mu^V$. This extension of the BD formula to the interacting field $W+Z$ is the subject of the next lemma.

Let $\mathcal{E}(Z, K)$ the symmetric bilinear form associated to $\mathcal{E}(Z)$, in such a way that

$$\mathcal{E}(Z+K) = \mathcal{E}(Z) + \mathcal{E}(K) + 2\mathcal{E}(Z, K).$$

**Lemma 13.** *Let V be admissible and denote by $\mu^V \in \bar{\mathcal{X}}$ a minimizer of $\check{F}^V$. Then we have the following direct variational description of $\mu^V$:*

$$\mathcal{W}^V(f) = -\log \mathbb{E}_{\mu^V}[e^{-f(W+Z)}] = \inf_{v \in \bar{\mathcal{Y}}(\mu^V)} \check{G}^{f,V}(v)$$

*with*

$$\check{G}^{f,V}(v) := \mathbb{E}_v[f(W+Z+K)] + \check{H}^V(v),$$

*and*

$$\check{H}^V(v) := \mathbb{E}_v[V(W+Z+K) - V(W+Z) + 2\mathcal{E}(Z,K) + \mathcal{E}(K)].$$

*Moreover $\nu$ is a minimizer of $\check{F}^V$ iff $\check{H}^V(v) \geqslant 0$ for all $v \in \bar{\mathcal{Y}}(\nu)$ and if V is differentiable in the directions of $H^1$ we also have the weak forward–backward SDE (wFBSDE):*

$$\mathbb{E}_v[\nabla V(W+Z) \cdot K + 2\mathcal{E}(Z, K)] = 0, \qquad v \in \bar{\mathcal{Y}}(\nu), \tag{14}$$

*and $\check{H}^V(v)$ takes the alternative form*

$$\check{H}^V(v) = \mathbb{E}_v[V(W+Z+K) - V(W+Z) - \nabla V(W+Z) \cdot K + \mathcal{E}(K)].$$

**Proof.** Recall that

$$\begin{aligned}
\mathcal{W}^V(f) &= \inf_{\mu \in \bar{\mathcal{X}}} \check{F}^{f+V} - \inf_{\mu \in \bar{\mathcal{X}}} \check{F}^V \\
&= \inf_{\mu \in \bar{\mathcal{X}}} \mathbb{E}_\mu[f(W+Z) + V(W+Z) + \mathcal{E}(Z)] - \mathbb{E}_{\mu^V}[V(W+Z) + \mathcal{E}(Z)] \\
&= \inf_{v \in \bar{\mathcal{Y}}(\mu^V)} \mathbb{E}_v[f(W+Z+K) + V(W+Z+K) + \mathcal{E}(Z+K)] \\
&\quad - \mathbb{E}_{\mu^V}[V(W+Z) + \mathcal{E}(Z)]
\end{aligned}$$

where in the last line we used Lemma 5 to justify the change of variables in the minimization problem. Note that for any $v \in \bar{\mathcal{Y}}(\mu^V)$ we have that

$$\mathbb{E}_{\mu^V}[V(W+Z) + \mathcal{E}(Z)] = \mathbb{E}_v[V(W+Z) + \mathcal{E}(Z)]$$

so

$$\mathcal{W}^V(f) = \inf_{v \in \bar{\mathcal{Y}}(\mu^V)} \left\{ \mathbb{E}_v[f(W+Z+K)] + \check{H}^V(v) \right\}.$$

The fact that $\check{H}^V(v) \geqslant 0$ iff $\mathrm{Law}_v(W, Z) = \nu$ is a minimizer of $\check{F}^V$ follows from its definition, indeed note that

$$\check{H}^V(v) = \check{F}^V(\mathrm{Law}_v(W, Z+K)) - \check{F}^V(\mathrm{Law}_v(W, Z)).$$

If V is differentiable in the direction of $H^1$ we have that for any minimizer $\nu$ of $\check{F}^V$ and any $v \in \bar{\mathcal{Y}}(\nu)$:

$$0 = \nabla \check{F}^V(v) = \mathbb{E}_v[\nabla V(W+Z) \cdot K + 2\mathcal{E}(Z, K)]$$

which allow also to simplify the expression of $\check{H}^V$. □

**Remark 14.** The interest of this formulation is two fold:

a) it describes a minimizer $\nu$ (our optimal control) with the first order condition (14) which we refer to as a weak forward–backward SDE which allows us to obtain useful informations as we remarked in the introduction.

b) it shows that the renormalized cost functional $\check{H}^V(v)$ can be made finite in the infinite volume by choosing localized perturbations $K$. This remark will be clearer by looking at how the determine the infinite volume limit of the exponential model in Section 4.

## 2 The quartic interaction

In this section we will consider the above general setting in the case where the interaction potential $V$ is given by

$$V_\Lambda(\phi) := \lambda \int_\Lambda [\![\phi^4]\!],$$

where $\Lambda \subseteq \mathbb{R}^2$ is a measurable bounded set (of sufficient regularity, see below) and $[\![\phi^n]\!]$ denotes the $n$-th Wick power of the Gaussian free field $\phi$ [17] and $\lambda > 0$ an arbitrary fixed constant. Let $\mathcal{W}_\Lambda := \mathcal{W}^{V_\Lambda}$ and $\theta_\Lambda := \theta^{V_\Lambda}$. The measure $\theta_\Lambda$ on $\mathcal{S}'(\mathbb{R}^2)$ is the $\varphi_2^4$ Euclidean quantum field with a finite volume cutoff.

We will prove that for a sufficiently large class of functions $f$ the family $(\mathcal{W}_\Lambda(f))_\Lambda$ is uniformly bounded as $\Lambda \nearrow \mathbb{R}^2$. This implies that the family of measures $(\theta_\Lambda)_\Lambda$ is tight on $\mathcal{S}'(\mathbb{R}^2)$. We will describe any accumulation point via a formula for its Laplace transform involving the wFBSDE obtained in the previous section.

**Remark 15.** Let us note that all the discussion in this section can be extended to general $P(\varphi)_2$ models with a positive coefficient in the highest order polynomial. Moreover also $O(N)$ invariant vector models could be handled in similar ways. Convexity of the classical potential (i.e. before Wick ordering) is only needed in Section 3 below to prove large deviations.

We will need the following enhancement $\mathbb{W}$ of the GFF $W$. Let

$$\mathbb{W}_t^2 := [\![W_t^2]\!], \quad \text{and} \quad \mathbb{W}_t^3 := [\![W_t^3]\!],$$

where here $[\![\cdot]\!]$ denotes the Wick powers with respect to the Gaussian field $W$. Then

$$\mathbb{W} := (W, \mathbb{W}^2, \mathbb{W}^3) \in L^p(\mathbb{P}, \mathfrak{S})$$

for all $p \geqslant 1$, where $\mathfrak{S} := C([0,1], (C^{-\delta}(\rho^{-1/4}))^3)$ for a fixed $\delta > 0$ which can be chosen as small as we like. By abuse of notation we will denote also $\mathbb{W} = (W, \mathbb{W}^2, \mathbb{W}^3)$ the canonical process on $\mathfrak{S}$.

In the following we will always consider a family of $(\Lambda)_{|\Lambda| \in \mathbb{R}}$ given by $\Lambda = |\Lambda|^{1/2} S$ for some fixed set $S$ of volume 1 with $C^1$-boundary. This allow to prove that $\mathbb{1}_\Lambda$ is sufficiently regular to make sense of quantities like $\int_\Lambda [\![\phi^4]\!]$, $\int_\Lambda [\![W_1^4]\!]$ and in general of products with Wick monomials of the GFF.

**Lemma 16.** *For any $p \geqslant 1$, $s < \frac{1}{p}$ and $\iota > \frac{1}{16}$ we have that $\mathbb{1}_\Lambda \in B_{p,p}^s(\rho^\iota)$, furthermore*

$$\sup_{|\Lambda|} \|\mathbb{1}_\Lambda\|_{B_{p,p}^s(\rho^\iota)} < \infty.$$

**Proof.** We have that $\nabla \mathbb{1}_\Lambda = \nu d\sigma$ where $d\sigma$ is the surface measure on $\Lambda$ and $\nu$ is the outer unit normal. It is known that $\|f\|_{B^s_{p,p}(\rho)} \lesssim \|f\|_{B^s_{p,p}(\rho)} + \|\nabla f\|_{B^{s-1}_{p,p}(\rho)}$ (see e.g [19], Lemma A.5). This implies that for any $s<1$, $\sup_{|\Lambda|} \|\mathbb{1}_\Lambda\|_{B^s_{1,1}(\rho^\iota)} < \infty$. Clearly $\sup_{|\Lambda|} \|\mathbb{1}_\Lambda\|_{L^\infty(\rho^\iota)} < \infty$ and by Besov space interpolation (see e.g. [19] Lemma A.3)

$$\|\mathbb{1}_\Lambda\|_{B^{\theta s}_{p,p}(\rho^\iota)} \lesssim \|\mathbb{1}_\Lambda\|^\theta_{B^s_{1,1}(\rho^\iota)} \|\mathbb{1}_\Lambda\|^{1-\theta}_{L^\infty(\rho^\iota)}$$

$\|\mathbb{1}_\Lambda\|_{L^\infty(\rho^\iota)}$ where

$$\frac{1}{p} = \theta. \qquad \square$$

**Lemma 17.** *The functional $\check{F}^\Lambda := \check{F}^{V_\Lambda}$ is admissible and for all $\mu \in \bar{\mathcal{X}}$*

$$\check{F}^\Lambda(\mu) = \check{F}^{\mathcal{V}_\Lambda}(\mu) = \mathbb{E}_\mu [\mathcal{V}_\Lambda(\mathbb{W}, Z) + \mathcal{E}(Z)],$$

*with*

$$\mathcal{V}_\Lambda(\mathbb{W}, Z) := 4\lambda \int_\Lambda \mathbb{W}_1^3 Z_1 + 6\lambda \int_\Lambda \mathbb{W}_1^2 (Z_1)^2 + 4\lambda \int_\Lambda W_1(Z_1)^3 + \lambda \int_\Lambda (Z_1)^4. \qquad (15)$$

**Proof.** It is known that

$$\mathbb{E}\big[e^{-q\lambda \int_\Lambda [\![W^4]\!]}\big] < \infty,$$

for any $q>0$ (see e.g [17] or [3]) and

$$\mathbb{E}[|V_\Lambda(W)|^p] \leqslant \mathbb{E}\Big[\Big|\|[\![W^4]\!]\|_{B^{-\delta}_{\infty,\infty}(\rho^\delta)}\Big|^p\Big] < \infty,$$

for some $\delta > 0$, which shows that $F$ is a tame functional. Now it remains to prove the equality in (15). For this it is enough to show that for $Z$ such that $\|Z\|_{H^1(\rho^{1/2})} < \infty$ $\mu$-almost surely, we have

$$\int_\Lambda [\![(W_1 + Z_1)^4]\!] = \int_\Lambda [\![W_1^4]\!] + 4\int_\Lambda [\![W_1^3]\!]Z_1 + 6\int_\Lambda [\![W_1^2]\!](Z_1)^2 + 4\int_\Lambda W_1(Z_1)^3 + \int_\Lambda (Z_1)^4.$$

To see this, recall that if $C^\infty_c \ni \psi_\varepsilon \to \delta_0$, is a standard mollifier, $W^\varepsilon = \psi_\varepsilon * W$ and $a_\varepsilon = 6\mathbb{E}[(W^\varepsilon_1(0))^2]$, $b_\varepsilon = 3\mathbb{E}[(W^\varepsilon_1(0))^2]^2$

$$\int_\Lambda [\![(W_1 + Z_1)^4]\!] = \lim_{\varepsilon \to 0} \int_\Lambda (\psi_\varepsilon * (W_1 + Z_1))^4 - a_\varepsilon \int_\Lambda (\psi_\varepsilon * (W_1 + Z_1))^2 - b_\varepsilon$$

$$= \lim_{\varepsilon \to 0} \int_\Lambda (W^\varepsilon_1)^4 - a_\varepsilon (W^\varepsilon_1)^2 - b_\varepsilon + 4\int_\Lambda \Big((W^\varepsilon_1)^3 - \frac{a_\varepsilon}{2}W^\varepsilon_1\Big)Z^\varepsilon_1$$

$$+ 6\int_\Lambda \Big((W^\varepsilon_1)^2 - \frac{a_\varepsilon}{6}\Big)(Z^\varepsilon_1)^2 + 4\int_\Lambda W^\varepsilon_1 (Z^\varepsilon_1)^3 + 4\int_\Lambda (Z^\varepsilon_1)^4$$

where we have denoted $Z^\varepsilon_1 = \psi_\varepsilon * Z_1$. Now it is well known that almost surely in $C^{-\delta}(\rho^\delta)$ we have $W^\varepsilon_1 \to W_1$ and

$$(W^\varepsilon_1)^4 - a_\varepsilon (W^\varepsilon_1)^2 - b_\varepsilon \to [\![W_1^4]\!], \quad \Big((W^\varepsilon_1)^3 - \frac{a_\varepsilon}{2}W^\varepsilon_1\Big) \to [\![W_1^3]\!], \quad \Big((W^\varepsilon_1)^2 - \frac{a_\varepsilon}{6}W^\varepsilon_1\Big) \to [\![W_1^2]\!].$$

Furthermore $Z^\varepsilon \to Z$ in $H^1$ which implies $(Z^\varepsilon_1)^2 \to Z_1^2, (Z^\varepsilon_1)^3 \to (Z^\varepsilon_1)^3$ in $W^{2\delta,1}$ for some $\delta > 0$, which gives the statement. From this calculation it also follows that $(W, Z) \to G(W+Z)$ is continuous on $\mathfrak{S} \times \mathfrak{H}_w$.

To check the lower semicontinuity of $\check{F}^{\mathcal{V}_\Lambda}$ recall from the Section 3 of [3] that

$$\mathcal{V}_\Lambda(\mathbb{W}, Z) + \mathcal{E}(Z) \geqslant -Q(\mathbb{W})$$

with $Q(\mathbb{W}) = C(\|\mathbb{W}\|_{\mathfrak{S}}^p + 1)$ for some $p \geqslant 1$. Furthermore observe that the function

$$(\mathbb{W}, Z) \mapsto \mathcal{V}_\Lambda(\mathbb{W}, Z) + \mathcal{E}(Z)$$

is lower semicontinuous on $\mathfrak{S} \times \mathfrak{H}_w$. Since the first marginal in $\bar{\mathcal{X}}$ is fixed, the Portmanteau theorem gives, for any $\mu_n \to \mu$ in $\bar{\mathcal{X}}$:

$$\begin{aligned}
\liminf_{n\to\infty} \check{F}^{\mathcal{V}_\Lambda}(\mu_n) &= \liminf_{n\to\infty} \mathbb{E}_{\mu_n}[\mathcal{V}_\Lambda(\mathbb{W}, Z) + \mathcal{E}(Z) + Q(\mathbb{W})] - \mathbb{E}_{\mu_n}[Q(\mathbb{W})] \\
&\geqslant \mathbb{E}_\mu[\mathcal{V}_\Lambda(\mathbb{W}, Z) + \mathcal{E}(Z) + Q(\mathbb{W})] - \mathbb{E}_\mu[Q(\mathbb{W})] \\
&= \check{F}^{\mathcal{V}_\Lambda}(\mu).
\end{aligned}$$

Finally, by a slight adaptation of the proof of Lemma 21 below we see that there exists $\delta > 0$ and $C_{\delta, |\Lambda|} < \infty$ such that

$$\check{F}^{\mathcal{V}_\Lambda}(\mu) \geqslant \delta \, \mathbb{E}_\mu[\mathcal{E}(Z) + \lambda \|\mathbb{1}_\Lambda Z_1\|_{L^4}^4] - C_{\delta, |\Lambda|},$$

for all $\mu \in \bar{\mathcal{X}}$, from which the strong coercivity condition (8) follows. $\square$

The admissibility of $\check{F}^\Lambda$ implies the variational representation for $\mathcal{W}_\Lambda$ and the existence of minimizers $\mu^{f,\Lambda}$ of $\check{F}^{f,\Lambda} := \check{F}^{f+\mathcal{V}_\Lambda}$ with the corresponding wFBSDE equation

**Lemma 18.** *For any bounded and continuous $f : \mathfrak{S}_0 \to \mathbb{R}$ there exists a minimizer $\mu^{f,\Lambda} \in \bar{\mathcal{X}}$ of $\check{F}^{f,\Lambda} := \check{F}^{f+\mathcal{V}_\Lambda}$. Moreover if $f$ is differentiable in the directions of $H^1$ any minimizer $\mu$ satisfies the wFBSDE equation*

$$\begin{aligned}
\mathbb{E}_\mu\left[4\lambda \int_\Lambda Z^3 K + 2\mathcal{E}(Z, K)\right] \\
= -\mathbb{E}_\mu\left[\nabla f(W + Z) \cdot K + \lambda \int_\Lambda (4\mathbb{W}^3 + 12\mathbb{W}^2 Z + 12\mathbb{W} Z^2) K\right]
\end{aligned} \tag{16}$$

*for any admissible control $K$.*

**Proof.** As we already remarked the existence of $\mu^{f,\Lambda}$ follows from the admissibility of $\check{F}^{f,\Lambda}$. The functional $Z \mapsto \mathcal{V}_\Lambda(\mathbb{W}, Z)$ is differentiable in $H^1$ and

$$\nabla \mathcal{V}_\Lambda(\mathbb{W}, Z) \cdot K = 4\lambda \int_\Lambda \mathbb{W}_1^3 K_1 + 12\lambda \int_\Lambda \mathbb{W}_1^2 Z_1 K_1 + 12\lambda \int_\Lambda W_1 (Z_1)^2 K_1 + 4\lambda \int_\Lambda (Z_1)^3 K_1.$$

Provided $Z \mapsto f(W + Z)$ is also differentiable, we have

$$0 = \nabla \check{F}^{f,\Lambda}(\mu)(K) = \mathbb{E}_\mu[\nabla f(W + Z) \cdot K + \nabla \mathcal{V}_\Lambda(\mathbb{W}, Z) \cdot K + 2\mathcal{E}(Z, K)]$$

which is (16). $\square$

**Remark 19.** Note that $Z \mapsto \mathcal{V}_\Lambda(\mathbb{W}, Z)$ is twice differentiable, with

$$\nabla^2 \mathcal{V}_\Lambda(\mathbb{W}, Z)(K, K') = 12\lambda \int_\Lambda \mathbb{W}_1^2 K_1' K_1 + 24\lambda \int_\Lambda W_1 Z_1 K_1' K_1 + 4\lambda \int_\Lambda (Z_1)^2 K_1 K_2.$$

Therefore if $f$ is twice differentiable then we have also the second order condition:

$$0 \leqslant \nabla^2 \check{F}^{f,\Lambda}(\mu)(K, K') = \mathbb{E}_\mu[\nabla^2 f(W + Z)(K, K') + \nabla^2 \mathcal{V}_\Lambda(\mathbb{W}, Z)(K, K') + 2\mathcal{E}(K', K)]$$

for any minimizer $\mu$ and any pair $K, K'$ of admissible controls.

We record also the normalized variational representation for $\mathcal{W}_\Lambda(f)$.

**Lemma 20.** *Denote by $\mu^\Lambda \in \bar{\mathcal{X}}$ a minimizer of $\check{F}^{0,\Lambda}$. Then we have*

$$\mathcal{W}_\Lambda(f) = \inf_{\upsilon \in \bar{\mathcal{Y}}(\mu^\Lambda)} \left\{ \mathbb{E}_\upsilon\left[\int_\Lambda f(W+Z+K)\right] + \check{H}^\Lambda(\rho) \right\}$$

*with*

$$\check{H}^\Lambda(\upsilon) := \mathbb{E}_\upsilon\left[\int_0^1 \int_{\mathbb{R}^2} \dot{Z}_s(m^2-\Delta)\dot{K}_s \mathrm{d}s + \frac{1}{2}\int_0^1 \int_{\mathbb{R}^2} \dot{K}_s(m^2-\Delta)\dot{K}_s \mathrm{d}s\right]$$

$$+ \mathbb{E}_\upsilon\left[4\lambda \int_\Lambda [\![(W+Z)^3]\!]\,K + 6\lambda \int_\Lambda [\![(W+Z)^2]\!]\,K^2 + 4\lambda \int_\Lambda [\![(W+Z)]\!]\,K^3 + \lambda \int_\Lambda K^4\right].$$

**Proof.** This follows from Lemma 13 by direct computation. □

Next we now show that the wFBSDE (16) gives us an a-priori bound uniform in $|\Lambda|$ on the minimizer in terms of weighted norms of $\mathbb{W}$. For this we first need some analytic bounds.

**Lemma 21.** *Assume that, for $\delta > 0$ small enough, either*

a) $\mathrm{supp}(\nabla f(W+Z)) \subset \Lambda$ and

$$\|\nabla f(W+Z)\|_{L^{4/3}(\rho^{3/4})} \leqslant \delta \|\mathbb{1}_\Lambda(W+Z)\|^3_{W^{-\varepsilon,4}(\rho^{1/4})} + C,$$

b)

$$\|\nabla f(W+Z)\|_{H^{-1}(\rho^{1/2})} \leqslant \delta \|W+Z\|_{H^{-\varepsilon}(\rho^{1/2})}.$$

*Then there exists $p \in (1,\infty)$ such that*

$$\mathbb{E}\left|\nabla f(W+Z)\cdot(Z\rho) + \lambda \int_\Lambda (\mathbb{W}^3\,\rho\,Z + 12\mathbb{W}^2\,\rho\,Z^2 + 12\,W\rho\,Z^3)\right|$$

$$\leqslant \delta\,\mathbb{E}\left(\lambda \int_\Lambda \rho\,Z^4 + \mathcal{E}(\rho^{1/2}Z)\right) + C.$$

**Proof.** First observe that we have that $\int_\Lambda \rho\,Z^4 = \|\mathbb{1}_\Lambda \rho^{1/4}\,Z^4\|^4_{L^4}$, and $\mathcal{E}(\rho^{1/2}Z) \geqslant C\|\rho^{1/2}\,Z\|^2_{H^1}$, for instance by fundamental theorem of calculus. By our assumptions on $f$ we have, in the first case

$$\begin{aligned}
|\nabla f(W+Z)\cdot(Z\rho)| &\leqslant \|\nabla f(W+Z)\|_{L^{4/3}(\rho^{3/4})}\|\mathbb{1}_\Lambda Z\|_{L^4(\rho^{1/4})} \\
&\leqslant \delta \|\mathbb{1}_\Lambda(W+Z)\|^3_{W^{-\varepsilon,4}(\rho^{1/4})}\|\mathbb{1}_\Lambda Z\|_{L^4(\rho^{1/4})} + C \\
&\leqslant C\delta\big(\|\mathbb{1}_\Lambda Z\|^4_{L^4(\rho^{1/4})} + \|\mathbb{1}_\Lambda W\|^3_{W^{-\varepsilon,4}(\rho^{1/4})}\|\mathbb{1}_\Lambda Z\|_{L^4(\rho^{1/4})}\big) + C \\
&\leqslant 2\delta\|\mathbb{1}_\Lambda Z\|^4_{L^4(\rho^{1/4})} + C\|\mathbb{1}_\Lambda W\|^4_{W^{-\varepsilon,4}(\rho^{1/4})} + C.
\end{aligned}$$

While in the second case

$$\begin{aligned}
|\nabla f(W+Z)\cdot(Z\rho)| &\leqslant \|\nabla f(W+Z)\|_{H^{-1}(\rho^{1/2})}\|Z\|_{H^1(\rho^{1/2})} \\
&\leqslant \delta\|W+Z\|_{H^{-\varepsilon}(\rho^{1/2})}\|Z\|_{H^1(\rho^{1/2})} \\
&\leqslant \delta\|W\|_{H^{-\varepsilon}(\rho^{1/2})}\|Z\|_{H^1(\rho^{1/2})} + \delta\|Z\|^2_{H^1(\rho^{1/2})} \\
&\leqslant \|W\|^2_{H^{-\varepsilon}(\rho^{1/2})} + 2\delta\|Z\|^2_{H^1(\rho^{1/2})}.
\end{aligned}$$

For the second term, by Cauchy–Schwarz we obtain

$$\left|\int_\Lambda \rho\,\mathbb{W}^3\,Z\right| \leqslant C_\delta\|\mathbb{1}_\Lambda \mathbb{W}^3\|^2_{H^{-1}(\rho^{1/2})} + \delta\|Z\|^2_{H^1(\rho^{1/2})}.$$

For the third term, we use Lemma A.7 from [19] to get

$$\left|\int_\Lambda \mathbb{W}^2 \rho Z^2\right| \leqslant \|\rho^{1/8}\mathbb{W}^2\|_{C^{-\varepsilon}} \|\rho^{7/8}\mathbb{1}_\Lambda Z^2\|_{B^{\varepsilon}_{1,1}}$$
$$\leqslant \|\rho^{1/8}\mathbb{W}^2\|_{C^{-\varepsilon}} \|\rho^{3/8}\mathbb{1}_\Lambda Z\|_{L^2} \|\rho^{1/2}Z\|_{H^{2\varepsilon}}$$
$$\leqslant C\|\rho^{1/8}\mathbb{W}^2\|_{C^{-\varepsilon}} \|\rho^{1/4}\mathbb{1}_\Lambda Z\|_{L^4} \|\rho^{1/2}Z\|_{H^{2\varepsilon}}$$
$$\leqslant C_\delta\|\rho^{1/8}\mathbb{W}^2\|^4_{C^{-\varepsilon}} + \delta(\|\rho^{1/4}\mathbb{1}_\Lambda Z\|^4_{L^4} + \|\rho^{1/2}Z\|^2_{H^{2\varepsilon}}).$$

For the final term we proceed similarly. By Lemma 16 we have $\sup_{|\Lambda|} \|\mathbb{1}_\Lambda\|_{B^s_{p,p}(\rho^\iota)} < \infty$ for any $s < 1/p$ and any $\iota > 1/16$. Then with $\varepsilon < 1/8$ and for any $\iota > 1/16$ we have by fractional Leibniz rule

$$\left|\int_\Lambda \rho W Z^3\right| \leqslant \|\rho^{1/8-\iota}W\|_{\mathcal{C}^{-\varepsilon}} \|\rho^{3/8+\iota}\mathbb{1}_\Lambda Z\|_{H^\varepsilon} \|\rho^{1/4}\mathbb{1}_\Lambda Z\|^2_{L^4}$$
$$\lesssim \|\rho^{1/8-\iota}W\|_{\mathcal{C}^{-\varepsilon}} \|\rho^{3/8}Z\|_{B^{1/2}_{8/3,8/3}} \|\rho^{1/4}\mathbb{1}_\Lambda Z\|^2_{L^4}$$
$$\leqslant C\|\rho^{1/8-\iota}W\|_{\mathcal{C}^{-\varepsilon}} \|\rho^{1/2}Z\|^{1/2}_{H^1} \|\rho^{1/4}\mathbb{1}_\Lambda Z\|^{5/2}_{L^4}$$
$$\leqslant C_\delta\|\rho^{1/8-\iota}W\|^8_{\mathcal{C}^{-\varepsilon}} + \delta(\|\rho^{1/2}Z\|^2_{H^1} + \|\rho^{1/4}\mathbb{1}_\Lambda Z\|^4_{L^4}).$$

and we can choose $\iota = 1/8 - \varepsilon > 1/16$ to get the claim. □

**Theorem 22.** *Assume that $f$ satisfies the conditions of Lemma 21. There exists a constant $C$ independent of $|\Lambda|$ such that, for any minimizer $\mu^{f,\Lambda} \in \bar{\mathcal{X}}$ of $\check{F}^{f,\Lambda}(\mu)$ we have*

$$\mathbb{E}_{\mu^{f,\Lambda}}\left[2\lambda\int_\Lambda \rho Z^4 + \mathcal{E}(\rho^{1/2}Z)\right] \leqslant C. \tag{17}$$

**Proof.** Note that the map $Z \mapsto \rho Z$ is continuous from $\mathfrak{H}$ to $\mathfrak{H}$. Take $K = \rho Z$ in (16) to get that any minimizer $\mu^{f,\Lambda} \in \bar{\mathcal{X}}$ satisfies

$$\mathbb{E}_{\mu^{f,\Lambda}}\left[4\lambda\int_\Lambda \rho Z^4 + \int_0^1\int_{\mathbb{R}^2}((m^2-\Delta)^{1/2}\rho^{1/2}\dot{Z}_s)^2 \mathrm{d}s\right]$$
$$= \mathbb{E}_{\mu^{f,\Lambda}}\left[\int_0^1\int_{\mathbb{R}^2}\dot{Z}_s[(m^2-\Delta),\rho^{1/2}]\rho^{1/2}\dot{Z}_s \mathrm{d}s\right] \tag{18}$$
$$+ \mathbb{E}_{\mu^{f,\Lambda}}\left[\nabla f(W+Z)(Z\rho) + \lambda\int_\Lambda \rho(4\mathbb{W}^3 Z + 12\mathbb{W}^2 Z^2 + 12\,WZ^3)\right].$$

From Lemma 23 and Lemma 21 below we deduce that the r.h.s. is controlled by a small fraction of the l.h.s. + a constant, from which we obtain the claim. □

**Lemma 23.**
$$\left|\int_0^1\int_{\mathbb{R}^2}\dot{Z}_s[(m^2-\Delta),\rho^{1/2}]\rho^{1/2}\dot{Z}_s \mathrm{d}s\right| \leqslant \varepsilon\left[\int_0^1\|\dot{Z}_s\|^2_{H^1(\rho^{1/2})}\mathrm{d}s\right]$$

**Proof.** We observe that $m^2$ obviously commutes with $\rho^{1/2}$, and

$$[\Delta, \rho^{1/2}] = (\Delta\rho^{1/2}) + \nabla\rho^{1/2}\cdot\nabla.$$

Using the assumptions (6) on the weight $\rho$ we deduce

$$\|[\Delta,\rho^{1/2}]f\|_{L^2} \leqslant \varepsilon(\|\rho^{1/2}f\|_{L^2} + \|\rho^{1/2}\nabla f\|_{L^2}). \qquad \square$$

Let us now show that these estimates give tightness of $(\theta_\Lambda)_\Lambda$ as $|\Lambda| \to \infty$ on a suitable weighted Sobolev space.

**Lemma 24.** *Consider $(\theta_\Lambda)_\Lambda$ as a family of measures on $H^{-\varepsilon}(\rho)$. We have for some $\delta > 0$*

$$\sup_\Lambda \left[ \int \exp(\delta \|\mathbb{1}_\Lambda \phi\|^4_{W^{-\varepsilon,4}(\rho)}) \mathrm{d}\theta_\Lambda + \int \exp(\delta \|\phi\|^2_{H^{-\varepsilon}(\rho)}) \mathrm{d}\theta_\Lambda \right] < \infty$$

*In particular $(\theta_\Lambda)_\Lambda$ is exponentially tight on $H^{-2\varepsilon}(\rho)$ and any accumulation point $\theta$ satisfies*

$$\int \exp(\delta \|\phi\|^4_{W^{-\varepsilon,4}(\rho)}) \mathrm{d}\theta < \infty.$$

**Proof.** Let $g_N \in C^\infty(\mathbb{R}, R)$ be such that $g_N(x) = x$ for $|x| < N$ and $g$. Define $f^N(\phi) = g_N\big(\delta \|\phi\|^2_{H^\varepsilon(\rho^{1/2})}\big)$ we have $\nabla f^N(\phi) = (g_N)'\big(\|\phi\|^2_{H^\varepsilon(\rho^{1/2})}\big) \delta \langle D \rangle^{-2\varepsilon} \phi$ so

$$\|\nabla f^N(W + Z)\|_{H^{-\varepsilon}} \leqslant \delta \|W + Z\|_{H^{-\varepsilon}},$$

and this bound is uniform in $N$. This implies by Lemma 21 that

$$\sup_\Lambda \mathbb{E}_{\mu^{f^N,\Lambda}}[\|Z\|^2_{H^{-\varepsilon}}] < \infty,$$

and

$$\log \int \exp(f^N) \mathrm{d}\theta_\Lambda \leqslant \mathbb{E}_{\mu^{0,\Lambda}}[f^N(W+Z)] + \mathbb{E}_{\mu^{f^N,\Lambda}}[f^N(W+Z)] \leqslant C < +\infty,$$

And taking sup over $N$ we get

$$\log \int \exp\big(\delta \|\phi\|^2_{H^\varepsilon(\rho^{1/2})}\big) \mathrm{d}\theta_\Lambda \leqslant C$$

a fact that implies exponential tightness of $(\theta_\Lambda)_\Lambda$. For the second statement we proceed similarly. Indeed if $\theta$ is an accumulation point of $(\theta_\Lambda)_\Lambda$ then

$$\int \exp(\delta \|\phi\|^4_{W^{-\varepsilon,4}(\rho)}) \mathrm{d}\theta \leqslant \sup_{\tilde\Lambda} \int \exp(\delta \|\mathbb{1}_{\tilde\Lambda} \phi\|^4_{W^{-\varepsilon,4}(\rho)}) \mathrm{d}\theta$$

$$\leqslant \sup_{\tilde\Lambda} \limsup_{\Lambda \to \infty} \int \exp(\delta \|\mathbb{1}_{\tilde\Lambda} \phi\|^4_{W^{-\varepsilon,4}(\rho)}) \mathrm{d}\theta_\Lambda$$

$$\leqslant \sup_{\tilde\Lambda} \sup_{\Lambda \supseteq \tilde\Lambda} \int \exp(\delta \|\mathbb{1}_{\tilde\Lambda} \phi\|^4_{W^{-\varepsilon,4}(\rho)}) \mathrm{d}\theta_\Lambda.$$

This gives the claim provided that for $f^N(\phi) = g^N\big(\delta \|\mathbb{1}_{\tilde\Lambda} \phi\|^4_{W^{-\varepsilon,4}(\rho^{1/4})}\big)$ we have

$$\sup_N \log \int \exp(f^N) \mathrm{d}\theta_\Lambda \leqslant C.$$

Indeed we compute $\nabla f^N(\phi) = (g^N)'\big(\|\mathbb{1}_{\tilde\Lambda} \phi\|^4_{W^{-\varepsilon,4}(\rho^{1/4})}\big) \delta \rho \mathbb{1}_\Lambda \langle D \rangle^{-\varepsilon}(\langle D \rangle^{-\varepsilon} \phi)^3$, from which we can check that

$$\|\nabla f^N(\phi)\|_{L^{4/3}(\rho^{3/4})} \leqslant C\delta \|\rho \mathbb{1}_\Lambda (\langle D \rangle^{-\varepsilon} \phi)^3\|_{L^{4/3}} \leqslant C\delta \|\rho^{1/4} \langle D \rangle^{-\varepsilon} \phi\|^3_{L^4} \leqslant C\delta \|\phi\|^3_{W^{-\varepsilon/2,4}},$$

and $\mathrm{supp}(\nabla f^N(\phi)) \subseteq \Lambda$. Applying Lemma 21 we get that the minimizer $\mu^{f^N,\Lambda}$ satisfies

$$\sup_{\tilde\Lambda} \sup_{\Lambda \supseteq \tilde\Lambda} \mathbb{E}_{\mu^{f^N,\Lambda}}\big[\|\mathbb{1}_{\tilde\Lambda} Z\|^4_{L^4(\rho^{1/4})}\big] \leqslant \sup_\Lambda \mathbb{E}_{\mu^{f^N,\Lambda}}\big[\|\mathbb{1}_\Lambda Z\|^4_{L^4(\rho^{1/4})}\big] =: C < \infty.$$

As a consequence, eq. (10) gives

$$\log \int \exp(f^N) \mathrm{d}\theta_\Lambda \leqslant \mathbb{E}_{\mu^{0,\Lambda}}[f^N(W+Z)] + \mathbb{E}_{\mu^{f^N,\Lambda}}[f^N(W+Z)] \leqslant 2C < +\infty,$$

which implies the claim. □

We have actually a more refined result: let $\mu^{f,\Lambda}$ be a minimizer for $\check{F}^{f,\Lambda}$. By Theorem 22 the family $(\mu^{f,\Lambda})_\Lambda$ is tight in $\mathcal{P}(\mathfrak{S} \times \hat{\mathfrak{H}}_w)$.

**Lemma 25.** *Denote by $\mu \in \mathcal{P}(\mathfrak{S} \times \hat{\mathfrak{H}}_w)$ any accumulation point of $(\mu^{f,\Lambda})_\Lambda$ as $|\Lambda| \to \infty$. Then the following equation holds:*

$$\mathbb{E}_\mu\bigg[\nabla f(W+Z) \cdot K + 4\lambda \int_{\mathbb{R}^2} [\![(W+Z)^3]\!]\, K + 2\mathcal{E}(Z,K)\bigg] = 0. \tag{19}$$

*for any $K: \mathfrak{S} \times \hat{\mathfrak{H}}_w \to \hat{\mathfrak{H}}_w$ continuous, adapted to the filtration generated by $(\mathbb{W}_t, Z_t)_{t\in[0,1]}$ and compactly supported.*

**Proof.** Let the support of $K$ be contained in $B_N(0)$ for some $N \geqslant 0$. Let $\chi \in C_c^\infty(\mathbb{R}^2)$ such that $\chi = 1$ on $B_N(0)$ so that $K = \chi K$. Then Lemma 21 implies that

$$\mathbb{E}_{\mu^{f,\Lambda}}[\|\chi Z\|_{H^1}^2 + \|\chi Z\|_{L^4}^4] \leqslant C,$$

where the constant depends on $N$ but not on $\Lambda$ and we have

$$\mathbb{E}_\mu\bigg[\nabla f(W+Z) \cdot K + 4\lambda \int_{\mathbb{R}^2} [\![(W+Z)^3]\!]\, K + 2\mathcal{E}(Z,K)\bigg]$$
$$= \mathbb{E}_\mu\bigg[\nabla f(W+Z) \cdot (\chi K) + 4\lambda \int_{\mathbb{R}^2} [\![(W+Z)^3]\!]\, \chi K + 2\mathcal{E}(Z,\chi K)\bigg]$$

Now observe that

$$4\lambda \int_{\mathbb{R}^2} [\![(W+Z)^3]\!]\, \chi K = H(\mathbb{W},Z,K) + 4\lambda \int Z^3 \chi K$$

where

$$H(\mathbb{W},Z,K) = 4\lambda \int_{\mathbb{R}^2} \mathbb{W}^3 \chi K + 12\lambda \int \mathbb{W}^2 Z \chi K + 12\lambda \int \mathbb{W} Z^2 \chi K.$$

Let $\tilde{\chi} \in C_c^\infty(\mathbb{R}^2)$ and $\tilde{\chi} = 1$ on the support of $\chi$. It is not hard to see from the definition that

$$|H(\mathbb{W},Z,K)| \lesssim \|\tilde{\chi}\mathbb{W}\|_\mathfrak{S}\big(1 + \|\chi^{1/3} Z\|_{L^4}^{4/3}\|\chi^{1/3} Z\|_{H^1}^{2/3}\big)\|\chi^{1/3} K\|_{B_{3,3}^{1/3}}$$

Indeed,

$$\bigg|\int_{\mathbb{R}^2} \mathbb{W}^3 \chi K\bigg| = \bigg|\int_{\mathbb{R}^2} \tilde{\chi}\mathbb{W}^3 \chi K\bigg| \lesssim \|\tilde{\chi}\mathbb{W}\|_\mathfrak{S} \|\chi K\|_{B_{3,3}^{1/3}},$$

while, by using interpolation from the first to the second line, we have

$$\bigg|\int_{\mathbb{R}^2} \mathbb{W}^2 Z \chi K\bigg| \lesssim \|\chi^{1/3}\mathbb{W}\|_\mathfrak{S} \|\chi^{1/3} Z\|_{B_{4/3,4/3}^{1/3}} \|\chi^{1/3} K\|_{B_{3,3}^{1/3}}$$
$$\lesssim \|\chi^{1/3}\mathbb{W}\|_\mathfrak{S} \|\chi^{1/3} Z\|_{L^4}^{2/3} \|\chi^{1/3} Z\|_{H^1}^{1/3} \|\chi^{1/3} K\|_{B_{3,3}^{1/3}}$$
$$\lesssim \|\tilde{\chi}\mathbb{W}\|_\mathfrak{S}\big(1 + \|\chi^{1/3} Z\|_{L^4}^{4/3}\|Z\|_{H^1}^{2/3}\big)\|\chi^{1/3} K\|_{B_{3,3}^{1/3}}$$

and similarly

$$\bigg|\int_{\mathbb{R}^2} W Z^2 \chi K\bigg|$$
$$\lesssim \|\tilde{\chi}\mathbb{W}\|_\mathfrak{S} \|\chi^{1/3} Z\|_{B_{4/3,4/3}^{1/3}}^2 \|\chi^{1/3} K\|_{B_{3,3}^{1/3}}$$
$$\lesssim \|\tilde{\chi}\mathbb{W}\|_\mathfrak{S} \|\chi^{1/3} Z\|_{L^4}^{4/3} \|\chi^{1/3} Z\|_{H^1}^{2/3} \|\chi^{1/3} K\|_{B_{3,3}^{1/3}}.$$

From this we get

$$\left|\nabla f(W+Z)\cdot(\chi K) + 4\lambda \int_{\mathbb{R}^2}[\![(W+Z)^3]\!]\,\chi\,K + 2\,\mathcal{E}(Z,\chi K)\right|$$
$$\lesssim \|\tilde{\chi}(W+Z)\|_{\mathcal{C}^{-\varepsilon}}\|\chi\,K\|_{H^1} + \left(\int_0^1 \|\chi^{1/2}\dot{Z}_t\|_{H^1}^2 dt\right)^{1/2}\left(\int_0^1 \|\chi^{1/2}\dot{K}_t\|_{H^1}^2 dt\right)^{1/2}$$
$$+ \|\tilde{\chi}\mathbb{W}\|_{\mathfrak{G}}\big(1 + \|\chi^{1/3}Z\|_{L^4}^{4/3}\|\chi^{1/3}Z\|_{H^1}^{2/3}\big)\|\chi^{1/3}K\|_{B_{3,3}^{1/3}}$$
$$+ \|\chi^{1/4}Z\|_{L^4}^3 \|\chi^{1/4}K\|_{L^4}.$$

The r.h.s. of this last equation is integrable by assumptions on $K$ and bounds on $Z$. From Lemmas 12 and 14 of [3] we know that there exist a family $(K^M)_{M \geqslant 1}$ adapted to the filtration generated by $(W_t)_t$ with

$$\begin{aligned}\mathbb{E}_\mu[\|\chi^{1/4}K\|_{L^4}^4 - \|\chi^{1/4}K^M\|_{L^4}^4] &\to 0, \\ \mathbb{E}_\mu\left[\left(\int_0^1 \|\chi^{1/2}(\dot{K}_t)\|_{H^1}^2 dt\right) - \int_0^1 \|\chi^{1/2}(\dot{K}_t^M)\|_{H^1}^2 dt\right] &\to 0,\end{aligned} \quad (20)$$

as $M \to \infty$. After possibly passing to a subsequence, we can assume that $\chi^{1/4}K^M \to \chi^{1/4}K$ in $L_w^4(\mu, L_w^4(\mathbb{R}^2))$ where $L_w^4$ denotes $L^4$ with the weak topology, and in $L_w^2(\mu, L^2(\mathbb{R}_+ \times \mathbb{R}^2))$. This together with the bound (20) implies that

$$\lim_{M\to\infty}\mathbb{E}_\mu\left[\nabla f(W+Z)\cdot(\chi K^M) + 4\lambda\int_{\mathbb{R}^2}[\![(W+Z)^3]\!]\,\chi\,K^M + 2\,\mathcal{E}(Z,\chi K^M)\right]$$
$$= \mathbb{E}_\mu\left[\nabla f(W+Z)\cdot(\chi K) + 4\lambda\int_{\mathbb{R}^2}[\![(W+Z)^3]\!]\,\chi\,K + 2\,\mathcal{E}(Z,\chi K)\right].$$

At this point we have reduced the problem to showing that

$$\lim_{\Lambda\to\infty}\mathbb{E}_{\mu^{f,\Lambda}}\left[\nabla f(W+Z)\cdot(\chi K^M) + 4\lambda\int_{\mathbb{R}^2}[\![(W+Z)^3]\!]\,\chi\,K^M + 2\,\mathcal{E}(Z,\chi K^M)\right]$$
$$= \mathbb{E}_\mu\left[\nabla f(W+Z)\cdot(\chi K^M) + 4\lambda\int_{\mathbb{R}^2}[\![(W+Z)^3]\!]\,\chi\,K^M + 2\,\mathcal{E}(Z,\chi K^M)\right]$$

since by assumption

$$\mathbb{E}_{\mu^{f,\Lambda}}\left[\nabla f(W+Z)\cdot(\chi K^M) + 4\lambda\int_{\mathbb{R}^2}[\![(W+Z)^3]\!]\,\chi\,K^M + 2\,\mathcal{E}(Z,\chi K^M)\right] = 0.$$

This follows from the weak convergence $\mu^{f,\Lambda} \to \mu$ and the bound

$$\mathbb{E}_{\mu^{f,\Lambda}}[\|\chi Z\|_{H^1}^2 + \|\chi Z\|_{L^4}^4] \leqslant C,$$

since we can approximate $Z$ by bounded admissible controls $(Z^L)_L$ defined as $Z_t^L = Z_{t \wedge T_L}$ with the stopping time

$$T_L := 1 \wedge \inf\{t \in [0,1] : \|\chi Z_t\|_{H^1}^2 + \|\chi Z_t\|_{L^4}^4 \geqslant L\},$$

in such a way that

$$\sup_{t \in [0,1]} \|\chi Z_t^L\|_{H^1}^2 + \|\chi Z_t^L\|_{L^4}^4 \leqslant L,$$

almost surely. Note that $T_L \to 1$ as $L \to \infty$ and therefore indeed $Z_t^L \to Z_t$ as $L \to \infty$. □

To summarize the result of this section: we obtained tightness of the finite volume $\varphi_2^4$ measures and described *any* possible infinite volume limit via a novel kind of stochastic equation, the wFBSDE (19). Note that it cannot be interpreted anymore as the first order condition of a functional. Unfortunately, due to the strongly non-convex character of the renormalized potential we have currently limited information on its solutions. One expects that in the "single phase region" [17] (e.g. where $m$ is large compared to $\lambda$) there should be uniqueness, and alternatively that multiple solutions should appear in "two phase region". Despite these limitations eq. (19) is useful and in the next section we show that it allows to prove large deviations for $\varphi_2^4$ in the semiclassical limit.

## 3 Large deviations

In this section we want to investigate the large deviations of $\varphi_2^4$ in the infinite volume in the semiclassical limit $\hbar \to 0$ where $\hbar > 0$ is a parameter which is inserted in the model as follows. For any $f : \mathcal{S}'(\mathbb{R}^2) \to \mathbb{R}$, let $\mathcal{W}_\Lambda^\hbar(f)$ be defined by

$$e^{-\frac{1}{\hbar}\mathcal{W}_\Lambda^\hbar(f)} := \int e^{-\frac{1}{\hbar}f(\phi)} \theta_\Lambda^\hbar(\mathrm{d}\phi),$$

where

$$\theta_\Lambda^\hbar(\mathrm{d}\phi) := \frac{\exp\left(-\frac{1}{\hbar}V_\Lambda(\phi)\right)\theta^\hbar(\mathrm{d}\phi)}{\int \exp\left(-\frac{1}{\hbar}V_\Lambda(\phi)\right)\theta^\hbar(\mathrm{d}\phi)}.$$

Here $\theta^\hbar$ is the Gaussian measure on $\mathcal{S}'(\mathbb{R}^2)$ with covariance $\hbar(m^2-\Delta)^{-1}$ and the Wick ordering is now taken with respect to $\theta^\hbar$. With some minor modifications to the derivation in the previous section we have the identity

$$\mathcal{W}_\Lambda^\hbar(f) = \inf_{\mu \in \bar{\mathcal{X}}} \check{F}^{f,\Lambda,\hbar}(\mu) - \inf_{\mu \in \bar{\mathcal{X}}} \check{F}^{0,\Lambda,\hbar}(\mu),$$

where

$$\check{F}^{f,\Lambda,\hbar}(\mu) := \mathbb{E}_{\mu \in \bar{\mathcal{X}}}[f(\hbar^{1/2}W + Z) + V_\Lambda(\hbar^{1/2}W + Z) + \mathcal{E}(Z)], \qquad \mu \in \bar{\mathcal{X}}.$$

Let $\mu^{f,\Lambda,\hbar}$ a minimizer for this admissible $\check{F}^{f,\Lambda,\hbar}$. Reasoning as in Lemma 24 above we have that the family $(\mu^{f,\Lambda,\hbar})_{\Lambda,\hbar}$ is exponentially tight for any bounded and continuous $f$. For each $\hbar > 0$ we pick a weakly convergent subsequence of $(\mu^{f,\Lambda,\hbar})_\Lambda$ with limit $\mu^{f,\hbar}$. We have then also $\theta_\Lambda^\hbar \to \theta^\hbar = \mathrm{Law}_{\mu^{f,\hbar}}(W_1 + Z_1)$ along the same subsequence and clearly, for all bounded continuous $f$ we have

$$\mathcal{W}_\Lambda^\hbar(f) \to \mathcal{W}^\hbar(f),$$

where $\mathcal{W}^\hbar(f)$ is now defined as

$$e^{-\frac{1}{\hbar}\mathcal{W}^\hbar(f)} := \int e^{-\frac{1}{\hbar}f(\phi)} \theta^\hbar(\mathrm{d}\phi).$$

In this setting where we do not have sufficient information to establish uniqueness of the limiting measures $(\mu^{f,\hbar})_\hbar$, nonetheless we are able to obtain a Laplace principle as $\hbar \to 0$ for the corresponding family $(\theta^\hbar)_\hbar$. Note that the family $(\theta^\hbar)_\hbar$ is exponentially tight by Lemma 24.

**Theorem 26.** *Any family of accumulation points $(\theta^\hbar)_\hbar$ of $(\theta_\Lambda^\hbar)_{\hbar,\Lambda}$ satisfies a Laplace principle on $\mathcal{S}'(\mathbb{R}^2)$ with rate function*

$$J(\phi) = \begin{cases} \lambda \int_{\mathbb{R}^2} \phi^4 + \frac{1}{2}\int_{\mathbb{R}^2} \phi(m^2-\Delta)\phi, & \text{for } \phi \in H^1(\mathbb{R}^2); \\ +\infty, & \text{otherwise.} \end{cases}$$

That is for $f\colon \mathcal{S}'(\mathbb{R}^2) \to \mathbb{R}$ *continuous and bounded*

$$\lim_{\hbar \to 0} \mathcal{W}^\hbar(f) = \inf_{\psi \in \mathcal{S}'(\mathbb{R}^2)} \{f(\psi) + J(\psi)\}.$$

**Proof.** Assume that $f$ is of the form

$$f(\psi) = C \int \rho((1-\Delta)^{-1/2}(\psi - \varphi))^2 \qquad (21)$$

for $C \geqslant 0$ and $\varphi \in C_c^\infty(\mathbb{R}^2)$.

It is enough to prove the theorem for this particular family of functions since we already know that $(\theta^\hbar)_\hbar$ are exponentially tight, and the family given by (21) is rate function determining in the sense of Definition 3.15 of [16], as it isolates points and is bounded below.

From (12) we have

$$\frac{\mathrm{d}}{\mathrm{d}\alpha} \mathcal{W}^\hbar_\Lambda(\alpha f) = \mathbb{E}_{\mu^{\alpha f, \Lambda, \hbar}}[f(\hbar^{1/2} W + Z)] = \frac{\int f(\phi) e^{-\frac{1}{\hbar} \alpha f(\phi)} \theta^\hbar_\Lambda(\mathrm{d}\phi)}{\int e^{-\frac{1}{\hbar} \alpha f(\phi)} \theta^\hbar_\Lambda(\mathrm{d}\phi)}.$$

By taking a further subsequence we can also have that $\mu^{\alpha f, \Lambda, \hbar} = \mu^{\alpha f, \Lambda, \hbar} \to \mu^{\alpha f, \hbar}$. We can then pass to the limit and obtain

$$\frac{\mathrm{d}}{\mathrm{d}\alpha} \mathcal{W}^\hbar(\alpha f) = \mathbb{E}_{\mu^{\alpha f, \hbar}}[f(\hbar^{1/2} W + Z)],$$

for any $f, \alpha$ and arbitrarily chosen accumulation point $\mu^{\alpha f, \hbar}$. Note that we do not have to do a diagonal argument over the uncountably many choices of $\alpha, f$: for any given $\alpha, f$ we need just to take a subsequence of the sequence chosen to have $(\theta^\hbar_\Lambda)_\Lambda$ converge in order to establish the equality.

Now recall that the minimizer $\mu^{f, \hbar}$ satisfies an infinite volume wFBSDE of the form

$$\begin{aligned} 0 &= \mathbb{E}_{\mu^{f,\hbar}} \bigg[ \nabla f(\hbar W + Z) \cdot H + \lambda \hbar^{3/2} \int \mathbb{W}^3 H + \lambda \hbar \int \mathbb{W}^2 Z H + \\ &\quad + \lambda \hbar^{1/2} \int W(Z)^2 H + \lambda \int Z^3 H + 2\, \mathcal{E}(Z, H) \bigg] \end{aligned} \qquad (22)$$

for any $H$ with compact support. Furthermore we can carry over the proof of Lemma 21, and obtain that

$$\sup_{\hbar < 1} \mathbb{E}_{\mu^{f,\hbar}} \big[ \|Z\|^4_{L^4(\rho^{1/4})} + \|Z\|^2_{H^1(\rho^{1/2})} \big] < \infty.$$

This implies that the family $(\mu^{f,\hbar})_\hbar$ is also tight and therefore we can consider any accumulation point $\mu^f$ as $\hbar \to 0$ which will satisfy the same estimate. By passing to the limit in (22) via by now standard estimates similar to those in Lemma 25 we obtain

$$0 = \mathbb{E}_{\mu^f} \bigg[ \nabla f(Z) \cdot H + \lambda \int Z^3 H + 2\, \mathcal{E}(Z, H) \bigg], \qquad (23)$$

valid for any adapted and compactly supported $H$. Note that this equation does not contain anymore the GFF. This equation has a unique solution. Indeed if $\nu$ is another solution we can couple the two solutions in $\rho \in \bar{\mathcal{Y}}(\mu^f)$ such that $\nu = \mathrm{Law}_\upsilon(W, Z + K)$ and have

$$0 = \mathbb{E}_\upsilon \bigg[ [\nabla f(Z + K) - \nabla f(Z)] \cdot H + \lambda \int [(Z+K)^3 - Z^3] H + 2\, \mathcal{E}(K, H) \bigg].$$

Now taking $H = \chi K$ for $\chi \in C_c^\infty(\mathbb{R}^2)$ and $\chi \geqslant 0$, we obtain

$$0 = \mathbb{E}_\nu\left[[\nabla f(Z+K) - \nabla f(Z)] \cdot (\chi K) + \lambda \int \chi[(Z+K)^3 - Z^3]K + 2\,\mathcal{E}(K, \chi K)\right].$$

By convexity of $f$ and $(\cdot)^4$ we have that

$$\mathbb{E}_\nu\left[[\nabla f(Z+K) - \nabla f(Z)] \cdot (\chi K) + \lambda \int \chi[(Z+K)^3 - Z^3]K\right] \geqslant 0.$$

Sending $\chi \to 1$ we obtain that

$$\mathbb{E}_\nu[\mathcal{E}(K,K)] \leqslant 0 \Rightarrow K = 0 \Rightarrow \nu = \mu^f.$$

Therefore $\mu^f$ is the unique minimizer of the admissible (convex) functional

$$\check{F}(\nu) := \mathbb{E}_\nu\left[f(Z) + \lambda \int Z^4 + \mathcal{E}(Z,Z)\right].$$

Let $\nu^\phi = \mathrm{Law}_\mathbb{P}(W, s \mapsto s\phi)$ for $\phi \in H^1$ then

$$\check{F}(\nu^\phi) = f(\phi) + \lambda \int \phi^4 + \frac{1}{2}\int \phi(m^2 - \Delta)\phi = f(\phi) + J(\phi).$$

Using the fact that the variational problem in the time variable for the functional

$$\mathcal{E}(Z,Z) = \frac{1}{2}\int_0^1 \int \dot{Z}_s (m^2 - \Delta) \dot{Z}_s \mathrm{d}s$$

is minimized for $\dot{Z}_s$ constant, we have for any $\nu$

$$\check{F}(\nu) = \mathbb{E}_\nu\left[f(Z) + \lambda \int Z^4 + \mathcal{E}(Z,Z)\right]$$

$$\geqslant \mathbb{E}_\nu\left[f(Z) + \lambda \int Z^4 + \frac{1}{2}\int Z_1(m^2 - \Delta) Z_1\right] = \mathbb{E}_\nu[f(Z) + J(Z)]$$

$$\geqslant \inf_{\phi \in H^1} [f(\phi) + J(\phi)].$$

As a consequence we obtain that

$$\inf_{\nu \in \bar{\mathcal{X}}} \check{F}(\nu) = \inf_{\phi \in H^1} [f(\phi) + J(\phi)],$$

and that $\mu^f = \nu^{\phi^f}$ for the unique $\phi^f$ which attain the convex deterministic minimization problem on the r.h.s..

To sum up, we have established that

$$\lim_{\hbar \to 0} \frac{\mathrm{d}}{\mathrm{d}\alpha} \mathcal{W}^\hbar(\alpha f) = \lim_{\hbar \to 0} \mathbb{E}_{\mu^{\alpha f, \hbar}}[f(\hbar^{1/2} W + Z)] = f(\phi^{\alpha f}). \tag{24}$$

Integrating eq. (24) and using dominated convergence

$$\lim_{\hbar \to 0} \mathcal{W}^\hbar(\alpha f) = \int_0^\alpha \lim_{\hbar \to 0} \frac{\mathrm{d}}{\mathrm{d}\beta} \mathcal{W}^\hbar(\beta f) \mathrm{d}\beta = \int_0^\alpha \lim_{\hbar \to 0} \mathbb{E}_{\mu^{\beta f, \hbar}}[f(\hbar^{1/2} W + Z)]\mathrm{d}\beta$$

$$= \int_0^\alpha f(\phi^{\beta f}) \mathrm{d}\beta =: \mathcal{W}^0(\alpha f).$$

This also implies that the derivative in $\alpha \in [0,1]$ of $\alpha \mapsto \mathcal{W}^0(\alpha f)$ exists and
$$\frac{\mathrm{d}}{\mathrm{d}\alpha}\mathcal{W}^0(\alpha f) = f(\phi^{\alpha f}).$$
To conclude it is enough to show that this coincides with
$$\frac{\mathrm{d}}{\mathrm{d}\alpha}\tilde{\mathcal{W}}^0(\alpha f),$$
where
$$\tilde{\mathcal{W}}^0(\alpha f) = \inf_{\psi \in H^1(\mathbb{R}^2)} \{\alpha f(\psi) + J(\psi)\}.$$
Recall that $\phi^{\alpha f}$ is the minimizer of $\phi \mapsto \alpha f(\phi) + J(\phi)$ so the Euler–Lagrange equation gives
$$\alpha \nabla f(\phi^{\alpha f}) \cdot h + \nabla J(\phi^{\alpha f}) \cdot h = 0.$$
At this point we note that
$$0 = \alpha[\nabla f(\phi^{(\alpha+\gamma)f}) - \nabla f(\phi^{\alpha f})] \cdot (\phi^{(\alpha+\gamma)f} - \phi^{\alpha f})$$
$$+ [\nabla J(\phi^{(\alpha+\gamma)f}) - \nabla J(\phi^{\alpha f})] \cdot (\phi^{(\alpha+\gamma)f} - \phi^{\alpha f})$$
$$+ \gamma \nabla f(\phi^{(\alpha+\gamma)f}) \cdot (\phi^{(\alpha+\gamma)f} - \phi^{\alpha f}).$$
By convexity of $f$ and strong convexity of $J$ we have
$$(\phi^{(\alpha+\gamma)f} - \phi^{\alpha f})^2$$
$$\lesssim \alpha[\nabla f(\phi^{(\alpha+\gamma)f}) - \nabla f(\phi^{\alpha f}) + \nabla J(\phi^{(\alpha+\gamma)f}) - \nabla J(\phi^{\alpha f})] \cdot (\phi^{(\alpha+\gamma)f} - \phi^{\alpha f}).$$
Now
$$\gamma \int \nabla f(\phi^{(\alpha+\gamma)f}) \cdot (\phi^{(\alpha+\gamma)f} - \phi^{\alpha f}) \mathrm{d}x \lesssim \gamma \|\nabla f(\phi^{(\alpha+\gamma)f})\|_{H^1} \|\phi^{(\alpha+\gamma)f} - \phi^{\alpha f}\|_{H^1}$$
so
$$\|(\phi^{(\alpha+\gamma)f} - \phi^{\alpha f})\|_{H^1}^2 \lesssim \gamma \|\nabla f(\phi^{(\alpha+\gamma)f})\|_{H^1} \|\phi^{(\alpha+\gamma)f} - \phi^{\alpha f}\|_{H^1}$$
so dividing by $\|\phi^{(\alpha+\gamma)f} - \phi^{\alpha f}\|_{H^1}$ we have
$$\|\phi^{(\alpha+\gamma)f} - \phi^{\alpha f}\|_{H^1} \leqslant \gamma \|\nabla f(\phi^{(\alpha+\gamma)f})\|_{H^1} \leqslant \gamma(\|\phi^{(\alpha+\gamma)f}\|_{H^1} + C) \lesssim \gamma.$$
We can then calculate
$$\tilde{\mathcal{W}}^0((\alpha+\gamma)f) - \tilde{\mathcal{W}}^0(\alpha f) = \gamma f(\phi^{(\alpha+\gamma)f}) + \alpha(f(\phi^{(\alpha+\gamma)f}) - f(\phi^{\alpha f}))$$
$$+ J(\phi^{(\alpha+\gamma)f}) - J(\phi^{\alpha f}).$$
Given that $\|\phi^{(\alpha+\gamma)f} - \phi^{\alpha f}\|_{H^1} \lesssim \gamma$ we have, since by Sobolev embedding $\|\cdot\|_{L^4}^4$ is locally Lipschitz in the $H^1$ norm
$$(\alpha(f(\phi^{(\alpha+\gamma)f}) - f(\phi^{\alpha f})) + J(\phi^{(\alpha+\gamma)f}) - J(\phi^{\alpha f}))$$
$$= \nabla f(\phi^{\alpha f})(\phi^{(\alpha+\gamma)f} - \phi^{\alpha f}) + \int (m^2 - \Delta)(\phi^{\alpha f})(\phi^{(\alpha+\gamma)f} - \phi^{\alpha f})$$
$$+ \lambda \int (\phi^{\alpha f})^3 (\phi^{(\alpha+\gamma)f} - \phi^{\alpha f}) + O(\gamma^2)$$
$$= O(\gamma^2),$$
where in the last line we used the Euler–Lagrange equation for $\phi^{\alpha f}$. This implies in turn
$$\frac{\mathrm{d}}{\mathrm{d}\alpha}\tilde{\mathcal{W}}^0(\alpha f) = \lim_{\gamma \to 0} \frac{1}{\gamma}(\tilde{\mathcal{W}}^0((\alpha+\gamma)f) - \tilde{\mathcal{W}}^0(\alpha f)) = \lim_{\gamma \to 0} \frac{1}{\gamma}\gamma f(\phi^{(\alpha+\gamma)f}) = f(\phi^{\alpha f})$$
which concludes the proof. $\square$

# 4 Exponential interaction

In this section we will discuss a variational description of the exponential interaction, that is the measure

$$\nu^{\xi,\exp}(\mathrm{d}\phi) = \frac{e^{-V^\xi(\phi)}\,\theta^0(\mathrm{d}\phi)}{\int e^{-V^\xi(\phi)}\,\theta^0(\mathrm{d}\phi)},$$

where now $V^\xi$ is given by

$$V^\xi(\phi) = \lambda \int_{\mathbb{R}^2} \xi [\![\exp(\beta\phi)]\!],$$

with some $\beta^2 \in (0, 8\pi)$ and where $M^\beta := [\![\exp(\beta\phi)]\!]$ is the Gaussian Multiplicative Chaos (GMC) and $\xi \in C_c^\infty(\mathbb{R}^2, \mathbb{R}_+)$ is a (smooth and compactly supported) spatial cutoff. We start with some preliminaries about the rigorous construction of the GMC.

**Proposition 27.** *Let $W$ be a GFF and $W_\varepsilon = \psi_\varepsilon * W$ where $\psi_\varepsilon \to \delta_0$ and $\psi_\varepsilon \in \mathcal{S}(\mathbb{R}^2)$ and let $\beta^2 \in (0, 8\pi)$. The Borel measures $(M^{\beta,\varepsilon})_\varepsilon$ on $\mathbb{R}^2$ given by*

$$M^{\beta,\varepsilon}(A) = \int_A \exp(-\beta W_\varepsilon(x) - \beta^2 \mathbb{E}[W_\varepsilon^2(x)])\mathrm{d}x, \qquad A \in \mathcal{B}(\mathbb{R}^2),$$

*converges almost surely set-wise to a measure*

$$M^\beta(A) = \lim_{\varepsilon \to 0} M^{\beta,\varepsilon}(A), \qquad A \in \mathcal{B}(\mathbb{R}^2). \tag{25}$$

*We call $M^\beta$ the Gaussian multiplicative chaos (GMC) with parameter $\beta$ associated to the GFF $W$. It enjoys the following properties*

a) *For $1 < p < 8\pi/\beta^2$ and for any compact set $A \subseteq \mathbb{R}^2$*

$$\sup_\varepsilon \mathbb{E}[|M^{\beta,\varepsilon}(A)|^p] < \infty.$$

*Consequently*

$$\mathbb{E}[|M^\beta(A)|^p] < \infty.$$

b) *For $1 < p < 8\pi/\beta^2$*

$$\mathbb{E}[|M^{\beta,\varepsilon}(B(0,\lambda))|^p]^{1/p} \leqslant C\lambda^{2-(p-1)\beta^2/4\pi} \mathbb{E}[|M^{\beta,\varepsilon}(B(0,1))|^p]^{1/p},$$

*where $B(x,r) \subseteq \mathbb{R}^2$ denotes the open ball centered in $x$ of radius $r$.*

**Proof.** For a proof of (25) see [7]. The property a) is due to Kahane [25], while the estimate in b) is due to Hoshino, Kawabi and Kusuoka [23]. □

From these properties we can determine easily the Besov regularity of $M^\beta$ following [23](Section 2).

**Proposition 28.** *Assume that $s > (\beta^2/(4\pi))(p-1)$, $1 < p < 8\pi/\beta^2$, then*

$$\mathbb{E}\left[ \|M^\beta\|^p_{B_{p,p}^{-s}(\rho^{1/2})} \right] < \infty.$$

**Proof.** We have

$$\begin{aligned}
\mathbb{E}\left[ \|M^{\beta,\varepsilon}\|^p_{B_{p,p}^{-s}(\rho^{1/2})} \right] &= \sum_{i \geqslant 0} 2^{-isp}\mathbb{E}\|\Delta_i M^{\beta,\varepsilon}\|^p_{L^p(\rho^{1/2})} \\
&\leqslant C\sum_{i \geqslant 0} 2^{-isp}\mathbb{E}|\Delta_i M^{\beta,\varepsilon}(0)|^p.
\end{aligned}$$

Recall that the Littlewood–Paley operators $(\Delta_i)_{i\geqslant -1}$ are defined by $\Delta_i = \psi(\mathrm{D}/2^i)$ for $i\geqslant 0$ and some Schwarz function $\psi\colon\mathbb{R}^2\to\mathbb{R}_{\geqslant 0}$. Now using translation invariance of the law of $M^{\beta,\varepsilon}$ and the fact that $\mathcal{F}^{-1}\psi$ is a Schwarz function we estimate

$$\begin{aligned}\mathbb{E}|\Delta_i M^{\beta,\varepsilon}(0)|^p &= 2^{2pi}\mathbb{E}\left|\int_{\mathbb{R}^2}(\mathcal{F}^{-1}\psi)(2^i x)M^{\beta,\varepsilon}(x)\mathrm{d}x\right|^p\\ &\leqslant \sum_{n\in\mathbb{Z}^2}(1+|n|^2)^{-2}\mathbb{E}\left|\int_{B(n,1)}M^{\beta,\varepsilon}(2^{-i}x)\mathrm{d}x\right|^p\\ &= \sum_{n\in\mathbb{Z}^2}(1+|n|^2)^{-2}2^{2pi}\mathbb{E}\left|\int_{B(n,2^{-i})}M^{\beta,\varepsilon}(x)\mathrm{d}x\right|^p\\ &\lesssim 2^{p(p-1)\beta^2 i/(4\pi)}.\end{aligned}$$

This gives the statement, since by weak convergence of $M^{\beta,\varepsilon}\to M^\beta$ we have

$$\mathbb{E}\|\Delta_i M^\beta\|_{L^p(\rho^{1/2})}^p \leqslant \limsup_{\varepsilon\to 0}\mathbb{E}\|\Delta_i M^{\beta,\varepsilon}\|_{L^p(\rho^{1/2})}^p.$$

$\square$

In this model we let $\mathfrak{S}:=\mathfrak{S}_0\times\mathcal{P}(\mathbb{R}^2)$ and $\mathbb{W}=(W,M^\beta)\colon\mathfrak{S}_0\to\mathfrak{S}$ will be the appropriate enhanced GFF. By the usual abuse of notation we will denote also $\mathbb{W}=(W,M^\beta)$ the canonical process on $\mathfrak{S}$.

We also need a slightly improved definition of the relaxed trial space $\bar{\mathcal{X}}^{\mathrm{exp}}$.

**Definition 29.** *Set*

$$\mathcal{X}^{\mathrm{exp}}:=\{\mu=\mathrm{Law}(W,M^\beta,Z)\in\mathcal{P}(\mathfrak{S}\times\hat{\mathfrak{H}}_w)\colon Z\in\hat{\mathfrak{H}}^a \text{ with } \mathbb{E}_\mu[\|Z\|_{\hat{\mathfrak{H}}}^2]<\infty\},$$

*and*

$$\bar{\mathcal{X}}^{\mathrm{exp}}:=\left\{\mu\in\mathcal{P}(\mathfrak{S}\times\hat{\mathfrak{H}}_w)\colon \exists(\mu_n)_n\subseteq\mathcal{X}\colon\mu_n\to\mu \text{ weakly}, \sup_{n\in\mathbb{N}}\mathbb{E}_{\mu_n}[\|Z\|_{\hat{\mathfrak{H}}}^2]<\infty\right\}$$

*We equip $\bar{\mathcal{X}}^{\mathrm{exp}}$ with the following topology: a sequence $(\mu_n)_n\subseteq\bar{\mathcal{X}}^{\mathrm{exp}}$ converges to $\mu$ if*

a) $\mu_n\to\mu$ *weakly in* $\mathcal{P}(\mathfrak{S}\times\hat{\mathfrak{H}}_w)$,

b) $\sup_{n\in\mathbb{N}}\mathbb{E}_{\mu_n}[\|Z\|_{\hat{\mathfrak{H}}}^2]<\infty.$

For technical reasons we first establish the variational formulation for the approximate potential $V^{\xi,\varepsilon}$ where

$$V^{\xi,\varepsilon}(\phi)=\lambda\int_{\mathbb{R}^2}\xi[\![\exp(\beta\psi_\varepsilon*\phi)]\!]$$

and the Wick ordering is now taken with respect to the covariance of $W_\varepsilon=\psi_\varepsilon*W$.

**Lemma 30.** *The potential $V^{\xi,\varepsilon}$ is admissible and*

$$\check{F}^{f,\xi,\varepsilon}(\mu):=\check{F}^{V^{\xi,\varepsilon}}(\mu)=\mathbb{E}_\mu\left[\int_{\mathbb{R}^2}\xi\exp(\beta(\psi_\varepsilon*Z))\mathrm{d}M^{\beta,\varepsilon}+\mathcal{E}(Z)\right]$$

*for all $\mu\in\bar{\mathcal{X}}^{\mathrm{exp}}$.*

**Proof.** Note that $M^{\beta,\varepsilon}$ is positive and that Proposition 28 implies that for $p$ close enough to 1 and $\frac{1}{p}+\frac{1}{q}=1$

$$\mathbb{E}[\exp(-qV^{\xi,\varepsilon})]+\mathbb{E}[|V^{\xi,\varepsilon}|^p]<\infty,$$

so $V^\xi$ is a tame functional. From the properties of the exponential it is clear that

$$\int_{\mathbb{R}^2} \xi[\![\exp(\beta(\psi_\varepsilon * Z) + \beta(\psi_\varepsilon * W_\varepsilon))]\!] \mathrm{d}x$$
$$= \int_{\mathbb{R}^2} \xi \exp(\beta(\psi_\varepsilon * Z))[\![\exp(\beta(\psi_\varepsilon * W_\varepsilon))]\!] \mathrm{d}x$$
$$= \int_{\mathbb{R}^2} \xi \exp(\beta(\psi_\varepsilon * Z)) \mathrm{d}M^{\beta,\varepsilon}.$$

We have for all $\mu \in \bar{\mathcal{X}}^{\mathrm{exp}}$:

$$\check{F}^{f,\xi,\varepsilon}(\mu) := \check{F}^{V^{\xi,f,\varepsilon}}(\mu)$$
$$= \mathbb{E}_\mu\Big[f(W+Z) + \int_{\mathbb{R}^2} \xi \exp(\beta(\psi_\varepsilon * Z)) \mathrm{d}M^{\beta,\varepsilon} + \mathcal{E}(Z)\Big]$$
$$\geqslant \mathbb{E}_\mu[\mathcal{E}(Z)] - C \geqslant \frac{1}{2}\mathbb{E}_\mu[\|Z\|_{\mathfrak{H}}^2] - C$$

from which also the coercivity of $\check{F}^{V^{\xi,f,\varepsilon}}$ is clear. Furthermore it is not hard to see that

$$(W, Z) \to \int_{\mathbb{R}^2} \xi[\![\exp(\beta(\psi_\varepsilon * Z) + \beta(\psi_\varepsilon * W))]\!] \mathrm{d}x$$

is continuous as a map $(\mathfrak{S}_0, H^1) \to \mathbb{R}$ (note hat here the renormalization if finite) and convergence of $W, Z$ in $(\mathfrak{S}_0, H^1)$ implies convergence in $L^\infty \times L^\infty$ of $\psi_\varepsilon * Z, \psi_\varepsilon * W$. To conclude admissibility it remains to show that $\check{F}^{V^{\xi,f,\varepsilon}}$ is lower semicontinuous on $\bar{\mathcal{X}}^{\mathrm{exp}}$. Since we have already shown the interaction to be continuous it remains to do this for $\mathcal{E}$. Indeed we have for any $\mu \in \bar{\mathcal{X}}^{\mathrm{exp}}$ a sequence $(\mu_n)_n$ such that $\mu_n \to \mu$ in $\bar{\mathcal{X}}^{\mathrm{exp}}$

$$\liminf_n \mathbb{E}_{\mu^n}[\mathcal{E}(Z)] \geqslant \mathbb{E}_\mu[\mathcal{E}(Z)],$$

this implies lower semicontinuity. $\square$

We now show that we can remove the ultraviolet cutoff in the variational formulation. We introduce the unregularized functional

$$\check{F}^{f,\xi}(\mu) = \check{F}^{V^{\xi,f}}(\mu) = \mathbb{E}_\mu\Big[f(W+Z) + \int_{\mathbb{R}^2} \xi \exp(\beta Z) \mathrm{d}M^\beta + \mathcal{E}(Z)\Big].$$

We will show that as $\varepsilon \to 0$ $\check{F}^{f,\xi,\varepsilon}$ $\Gamma$-converges to $\check{F}^{f,\xi}$. We briefly recall the notion of $\Gamma$-convergence, see e.g [10] for a nice introduction.

**Definition 31.** *Let $\mathcal{T}$ be a topological space and let $F, F_n: \mathcal{T} \to (-\infty, \infty]$. We say that the sequence of functionals $(F_n)_n$ $\Gamma$-converges to $F$ iff*

   i. *For every sequence $x_n \to x$ in $\mathcal{T}$*

$$F(x) \leqslant \liminf_{n \to \infty} F_n(x_n);$$

   ii. *For every point $x$ there exists a sequence $x_n \to x$ (called a recovery sequence) such that*

$$F(x) \geqslant \limsup_{n \to \infty} F_n(x_n).$$

**Lemma 32.** *$\check{F}^{f,\xi,\varepsilon}$ $\Gamma$-converges to $\check{F}^{f,\xi}$ on $\bar{\mathcal{X}}$ as $\varepsilon \to 0$.*

**Proof.** We first show the liminf inequality. We have for sequence $(\mu_\varepsilon)_\varepsilon$ such that $\mu_\varepsilon \to \mu$ in $\bar{\mathcal{X}}^{\mathrm{exp}}$

$$\liminf_{\varepsilon \to 0} \mathbb{E}_{\mu^\varepsilon}[\mathcal{E}(Z)] \geqslant \mathbb{E}_\mu[\mathcal{E}(Z)], \qquad \lim_{\varepsilon \to 0} \mathbb{E}_{\mu^\varepsilon}[f(W+Z)] = \mathbb{E}_\mu[f(W+Z)].$$

Now let $g^N: \mathbb{R} \to \mathbb{R}$ be defined by $g^N(x) = x$ if $x < N$, $g^N(x) = N$ if $x > N$. Then $g^N \leqslant N$, $g^N$ is Lipschitz with constant 1 and $g^N(x) \to x$. Set $Z_t^N = g^N(Z_t)$ and $\mu_{N,\varepsilon} = \mathrm{Law}_{\mu^\varepsilon}(W, Z^N)$. Then by monotone convergence

$$\lim_{N \to \infty} \mathbb{E}_\mu \int \xi \exp(\beta Z^N) \mathrm{d} M^\beta$$
$$= \mathbb{E}_\mu \int \xi \exp(\beta Z) \mathrm{d} M^\beta$$

It remains to prove that

$$\lim_{\varepsilon \to 0} \mathbb{E}_{\mu^\varepsilon} \int \xi \exp(\beta Z^{N,\varepsilon}) \mathrm{d} M^{\beta,\varepsilon} = \mathbb{E}_\mu \int \xi \exp(\beta Z^N) \mathrm{d} M^\beta,$$

then the statement will follow from

$$\mathbb{E}_\mu \int \xi \exp(\beta Z) \mathrm{d} M^\beta = \lim_{N \to \infty} \mathbb{E}_\mu \int \xi \exp(\beta Z^N) \mathrm{d} M^\beta$$
$$= \lim_{N \to \infty} \lim_{\varepsilon \to \infty} \mathbb{E}_{\mu^\varepsilon} \int \xi \exp(\beta Z^{N,\varepsilon}) \mathrm{d} M^{\beta,\varepsilon}$$
$$\leqslant \lim_{\varepsilon \to 0} \mathbb{E}_{\mu^\varepsilon} \int \xi \exp(\beta Z^\varepsilon) \mathrm{d} M^{\beta,\varepsilon}.$$

We have that $\|\exp(\beta Z^N)\|_{L^\infty} \leqslant \exp(\beta N)$ and by chain rule

$$\|\rho \nabla \exp(\beta Z^N)\|_{L^2} = \|\rho \beta \nabla Z^N \exp(\beta Z^N)\|_{L^2} \lesssim \exp(\beta N) \|\rho \nabla Z^N\|_{L^2}.$$

By interpolation

$$\|\exp(\beta Z^{N,\varepsilon})\|_{B_{p',p'}^\delta(\rho)} \lesssim \|\exp(\beta Z^{N,\varepsilon})\|_{L^\infty}^{(1-2\delta)} \|\exp(\beta Z^{N,\varepsilon})\|_{H^1(\rho)}^{2\delta},$$

provided $2/p' > \delta$. Now let $p = p'/(p'-1)$. Note that we can choose $\delta$ such that $\frac{8\pi}{\beta^2} > p$ and $\frac{\beta^2}{4\pi}(p-1) < \delta$ and $2/p' < \delta$ by taking

$$\frac{2}{p}(p-1) > \delta > \frac{\beta^2}{4\pi}(p-1).$$

Since $M^{\beta,\varepsilon} \to M^\beta$ in $B_{p,p}^s$ by Proposition 28 this establishes that $\int \xi \exp(\beta Z^{N,\varepsilon}) \mathrm{d} M^{\beta,\varepsilon} \to \int \xi \exp(\beta Z^N) \mathrm{d} M^\beta$ almost surely. Since

$$\int \xi \exp(\beta Z^{N,\varepsilon}) \mathrm{d} M^{\beta,\varepsilon} \lesssim e^{\beta N} \int \xi \mathrm{d} M^{\beta,\varepsilon},$$

and the r.h.s is uniformly integrable in the probability variable, we have by dominated convergence

$$\lim_{\varepsilon \to 0} \mathbb{E}_{\mu^\varepsilon} \int \xi \exp(\beta Z^{N,\varepsilon}) \mathrm{d} M^{\beta,\varepsilon} = \mathbb{E}_\mu \int \xi \exp(\beta Z^N) \mathrm{d} M^\beta,$$

and we can conclude. To show the limsup inequality, for given $\mu$ we need to construct $\mu^\varepsilon$ such that

$$\limsup_{\varepsilon \to 0} F^{V^{\xi,f},\varepsilon}(\mu^\varepsilon) \leqslant F^{V^{\xi,f}}(\mu).$$

To do this we proceed similarly as in the first step. Define $Z^{N,\delta} = g^N(\psi_\delta * Z)$, $\mu^{N,\delta} = \mathrm{Law}(W, M^\beta, Z^{N,\delta})$, where $\psi_\delta$ is a sequence of standard mollifiers. Then

$$M^{\beta,\varepsilon} \to M^\beta$$

weakly so
$$\mathbb{E}\int \exp(\beta Z^{N,\delta})\mathrm{d}M^{\beta,\varepsilon} \to \mathbb{E}\int \exp(\beta Z^{N,\delta})\mathrm{d}M^{\beta}$$
since $Z^{N,\delta}$ is bounded and continuous. Furthermore $\mathbb{E}_{\mu^\varepsilon}[\mathcal{E}(Z^{N,\delta})] \to \mathbb{E}_\mu[\mathcal{E}(Z^{N,\delta})]$. This implies that
$$\lim_{\varepsilon\to 0} F^{V^{\xi,\varepsilon}}(\mu^{N,\delta}) = F^{V^\xi}(\mu^{N,\delta})$$
Now by dominated/monotone convergence we have
$$\lim_{N\to\infty}\lim_{\delta\to 0}\mathbb{E}_\mu \int \exp(\beta Z^{N,\delta})\mathrm{d}M^\beta = \mathbb{E}_\mu \int \exp(\beta Z)\mathrm{d}M^\beta,$$
and
$$\lim_{N\to\infty}\lim_{\delta\to 0}\mathbb{E}_\mu[\mathcal{E}(Z^{N,\delta})] = \mathbb{E}_\mu[\mathcal{E}(Z)] \qquad \lim_{N\to\infty}\lim_{\delta\to 0}\mathbb{E}_\mu[f(W+Z^{N,\delta})] = \mathbb{E}_\mu[f(W+Z)],$$
so taking a suitable diagonal sequence implies the statement. $\square$

**Corollary 33.**
$$\lim_{\varepsilon\to 0}\inf_{\mu\in\bar{\mathcal{X}}} \check{F}^{V^{\xi,\varepsilon}} = \inf_{\mu\in\bar{\mathcal{X}}}\check{F}^{V^\xi}(\mu)$$
*and the infimum on the r.h.s is attained.*

**Proof.** The statement will follow from the fundamental theorem of $\Gamma$-convergence (see e.g. [14]), once we have established that $\check{F}^{V^{\xi,\varepsilon}}$ is equicoercive on $\bar{\mathcal{X}}^{\exp}$, but this follows from
$$\check{F}^{V^{\xi,f,\varepsilon}}(\mu) \geqslant \mathbb{E}_\mu[\mathcal{E}(Z)] - C. \qquad \square$$

We can now relate this back to our general framework:

**Lemma 34.** *The functional $\check{F}^{f,\xi}$ has a minimizer $\mu^{f,\xi}\in\bar{\mathcal{X}}^{\exp}$ and*
$$\mathcal{W}^{\xi,\exp}(f) := -\log\int \exp(-f(\phi))\mathrm{d}\nu^\xi = \min_{\mu\in\bar{\mathcal{X}}^{\exp}}\check{F}^{f,\xi}(\mu) = \check{F}^{f,\xi}(\mu^{f,\xi}).$$

*Furthermore for any minimizer $\mu\in\bar{\mathcal{X}}^{\exp}$ the following wFBSDE*
$$\mathbb{E}_\mu\left[\lambda\beta\int \xi\exp(\beta Z)\,K\mathrm{d}M^\beta + 2\mathcal{E}(Z,K)\right] = \mathbb{E}_\mu\left[\int \nabla f(W+Z)\,K\right]$$
*holds for any admissible control $K$.*

**Proof.** This follows *mutatis mutandis* as in Lemma 18. $\square$

**Lemma 35.** *The functional $\check{F}^{0,\xi}$ has a unique minimizer $\mu^{0,\xi}$ in $\bar{\mathcal{X}}$. Furthermore*
$$Z \leqslant 0,$$
*$\mu^{0,\xi}$ almost surely, and*
$$\sup_\xi \mathbb{E}_{\mu^{0,\xi}}[\|Z\|_{\mathfrak{H}}^2] < \infty.$$

**Proof.** Observe that $F^{0,\xi}(Z)$ is convex in $Z$. Now let $\mu,\nu$ be two minimizers for $\check{F}^{0,\xi}$. Recall from Lemma 5 that there exists $v\in\bar{\mathcal{Y}}(\mu)$ such that $\mathrm{Law}_v(W,Z+K)=\nu$. Now
$$\mathbb{E}_v\left[\lambda\int \xi\exp(\beta(Z+K))\mathrm{d}M^\beta + \mathcal{E}(Z+K)\right] - \mathbb{E}_v\left[\lambda\int \xi\exp(\beta Z)\mathrm{d}M^\beta + \mathcal{E}(Z)\right]$$
$$\geqslant \mathbb{E}_v\left[\lambda\beta\int \xi\exp(\beta Z)K\mathrm{d}M^\beta + 2\mathcal{E}(Z,K) + \mathcal{E}(K)\right].$$

By Lemma 34 we have that
$$\mathbb{E}_v\left[\lambda\beta\int \xi\exp(\beta Z)K\mathrm{d}M^\beta + 2\mathcal{E}(Z,K)\right]=0,$$
so
$$\mathbb{E}[\mathcal{E}(K)] \leqslant \check{F}^{f,\xi}(\mu) - \check{F}^{f,\xi}(\nu) = 0,$$
which implies $K=0$ and therefore $\nu=\mu$. For the second statement observe that the mapping $A\colon Z\to -|Z|$ leaves the energy $\mathcal{E}(Z)$ invariant by chain rule, while
$$\int \xi\exp(\beta A(Z))\mathrm{d}M^\beta \leqslant \int \xi\exp(\beta Z)\mathrm{d}M^\beta,$$
so if $\mu^{0,\xi}$ is a minimizer, so is $\mathrm{Law}_{\mu^{0,\xi}}(W,A(Z))$. By uniqueness we can conclude. For the third claim we test the wFBSDE with $K=\rho Z$ and obtain
$$\mathbb{E}_\mu\left[4\lambda\beta\int \xi\exp(\beta Z)\rho Z\mathrm{d}M^\beta + 2\mathcal{E}(Z,\rho Z)\right] = \mathbb{E}_\mu\left[\int \nabla f(W+Z)\rho Z\right].$$
We can again estimate as in the proof of Lemma 21
$$2\mathcal{E}(Z,\rho Z) \geqslant \frac{1}{2}\|Z\|_{\hat{\mathfrak{H}}}^2,$$
and
$$\int \nabla f(W+Z)\rho Z \lesssim \|Z\|_{\hat{\mathfrak{H}}}.$$
Furthermore observe that
$$\int_{\mathbb{R}^2} \xi\exp(\beta Z)\rho Z\mathrm{d}M^\beta \geqslant -C\int_{\mathbb{R}^2} \xi\rho\,\mathrm{d}M^\beta,$$
since $\exp(\beta Z)Z$ is bounded below. Summing up, we have the estimate
$$\frac{1}{2}\mathbb{E}_{\mu^{0,\xi}}[\|Z\|_{\hat{\mathfrak{H}}}^2] \lesssim \mathbb{E}[\|Z\|_{\hat{\mathfrak{H}}}^2]^{1/2} + \mathbb{E}_\mu\left[\int_{\mathbb{R}^2} \xi\rho\,\mathrm{d}M^\beta\right] \lesssim \mathbb{E}[\|Z\|_{\hat{\mathfrak{H}}}^2]^{1/2} + 1$$
where in the last bound we used Lemma 28. $\square$

**Lemma 36.** *Let $\mu^\xi \in \bar{\mathcal{X}}^{\exp}$. For $v \in \bar{\mathcal{Y}}(\mu^\xi)$ we have*
$$\check{H}^{\xi,\exp}(v) := \check{H}^{V^\xi}(v) = \mathbb{E}_v[2\mathcal{E}(K,Z) + \mathcal{E}(K)]$$
$$+ \lambda\mathbb{E}_v\left[\int \xi\exp(\beta(Z+K))\mathrm{d}M^\beta - \int \xi\exp(\beta Z)\mathrm{d}M^\beta\right].$$
*Then $\check{H}^{\xi,\exp}(v) \geqslant 0$ for all $v \in \bar{\mathcal{Y}}(\mu^\xi)$ iff $\mu^\xi$ is a minimizer of $\check{F}^{0,\xi}$. Moreover in this case, thanks to the wFBSDE for $Z$ we have*
$$\check{H}^{\xi,\exp}(v) = \mathbb{E}_v\left[\mathcal{E}(K) + \lambda\int \xi\exp(\beta Z)(\exp(\beta K) - 1 - \beta K)\mathrm{d}M^\beta\right],$$
*and*
$$\mathcal{W}_\xi^{\exp}(f) = \inf_{v\in\bar{\mathcal{Y}}(\mu^\xi)} \check{G}^{f,\xi,\exp}(v), \quad \check{G}^{f,\xi,\exp}(v) := \mathbb{E}_v[f(W+Z+K)] + \check{H}^{\xi,\exp}(v).$$

*Any minimizer $v \in \bar{\mathcal{Y}}(\mu^\xi)$ of $\check{G}^{f,\xi,\exp}$ satisfies the perturbed wFBSDE*
$$\begin{aligned}\mathbb{E}_v\bigg[\mathcal{E}(K) + \lambda\beta\int \xi\exp(\beta Z)(\exp(\beta K) - 1)H\mathrm{d}M^\beta\bigg] \\ = \mathbb{E}_v[\nabla f(W+Z+K)\cdot H].\end{aligned} \quad (26)$$

The proof of Lemma 36 is completely analogous to Lemma 20 and Lemma 18, thus we omit the details. Thanks to the convexity of the interaction we have a good control of the locality properties of the solutions to the wFBSDE (26).

**Lemma 37.** *Assume that $N > 0$ and*
$$\sup_{\varphi \in L^2(\mathbb{R}^2)} \|\nabla f(\varphi)\|_{H^{-1}(\rho^{-N})} < \infty.$$

*Then for any minimizer $\upsilon^f$ of $\check{G}^{f,\xi,\exp}(\rho)$ we have*
$$\mathbb{E}_{\upsilon^f}[\|K\|_{H^1(\rho^{-N})}^2] < \infty.$$

**Proof.** Testing equation (26) with $H = \rho^{-2N}K$ we have

$$\mathbb{E}_\upsilon\left[\frac{1}{2}\int_0^1 \int_{\mathbb{R}^2} \rho^{-2N}\dot{K}_s(m^2-\Delta)\dot{K}_s \mathrm{d}s\right]$$
$$+\lambda\beta\mathbb{E}_\upsilon\left[\int \xi\exp(\beta Z)(\exp(\beta K)-1)\rho^{-2N}K\mathrm{d}M^\beta\right]$$
$$= \mathbb{E}_\upsilon\left[\int \rho^{-2N}\nabla f(W+Z+K)K\right] \leqslant C\mathbb{E}_\upsilon[\|K\|_{H^1(\rho^{-N})}^2]^{1/2}$$

Now

$$\mathbb{E}_\upsilon\left[\frac{1}{2}\int_0^1 \int_{\mathbb{R}^2} \rho^{-2N}\dot{K}_s(m^2-\Delta)\dot{K}_s \mathrm{d}s\right]$$
$$= \mathbb{E}_\upsilon\left[\frac{1}{2}\int_0^1 \int_{\mathbb{R}^2} \rho^{-2N}((m^2-\Delta)^{-1/2}\dot{K}_s)^2 \mathrm{d}s\right]$$
$$+\mathbb{E}_\upsilon\left[\frac{1}{2}\int_0^1 \int_{\mathbb{R}^2}[(m^2-\Delta)^{1/2},\rho^{-2N}]\dot{K}_s(m^2-\Delta)^{1/2}\dot{K}_s \mathrm{d}s\right]$$

and we can compute, as in the proof of Lemma 21,

$$[(m^2-\Delta)^{1/2},\rho^{-2N}]\dot{K}_s \leqslant \delta\rho^{-2N}\dot{K}_s.$$

Furthermore observe that

$$\mathbb{E}_\upsilon\left[\int \xi\exp(\beta Z)(\exp(\beta K)-1)\rho^{-2N}K\mathrm{d}M^\beta\right] \geqslant 0,$$

Putting everything together we have

$$(1/2-\delta)\mathbb{E}_\upsilon\|K\|_{H^1(\rho^{-N})}^2 = (1/2-\delta)\mathbb{E}_\upsilon\left[\frac{1}{2}\int_0^1 \int_{\mathbb{R}^2}\rho^{-2N}((m^2-\Delta)^{-1/2}\dot{K}_s)^2\mathrm{d}s\right]$$
$$\leqslant C\mathbb{E}_\upsilon\|K\|_{H^1(\rho^{-N})}$$

which implies the statement. □

**Definition 38.** *For $N > 0$ large (to be chosen suitably below) denote by $\mathcal{Y}^M(\mu)$ the set of measures $\upsilon \in \mathcal{P}(\mathfrak{S} \times \hat{\mathfrak{H}}_w \times \mathfrak{H}_w(\rho^{-N}))$ such that*

  a) $\mathrm{Law}_\upsilon(W, M^\beta, Z) = \mu$

  b) $\exists (Z^n)_n, (K^n)_n \subseteq \hat{\mathfrak{H}}^a$ *for which*
$$\mathrm{Law}(W, M^{\beta,n}, Z^n, K^n) \to \upsilon \text{ and } \sup_{n \in \mathbb{N}}\mathbb{E}[\|Z^n\|_{\hat{\mathfrak{H}}}^2 + \|K^n\|_{\mathfrak{H}(\rho^{-N})}^2] \leqslant 2M.$$

Let $\mathcal{Y}^M := \bigcup_{\mu \in \bar{\mathcal{X}}^{\exp}} \mathcal{Y}^M(\mu)$. We say that $(\upsilon^n)_n \subseteq \mathcal{Y}^M$ converges to $\upsilon \in \mathcal{Y}^M$ if $\upsilon^n \to \upsilon$ weakly in $\mathcal{P}(\mathfrak{S} \times \hat{\mathfrak{H}}_w \times \mathfrak{H}_w(\rho^{-N}))$ and

$$\sup_{n \in \mathbb{N}} \mathbb{E}_{\upsilon^n}[\|Z\|^2_{\hat{\mathfrak{H}}_w} + \|K\|^2_{\mathfrak{H}(\rho^{-N})}] \leqslant 2M.$$

**Remark 39.** Note that any sequence in $\mathcal{Y}^M$ has a convergent subsequence.

**Corollary 40.** *There exists a sufficiently large $M = M(f)$ (not depending on $\xi$) such that*

$$\inf_{\upsilon \in \mathcal{Y}(\mu^{0,\xi})} \check{G}^{f,\xi,\exp}(\upsilon) = \inf_{\upsilon \in \mathcal{Y}^M(\mu^{0,\xi})} \check{G}^{f,\xi,\exp}(\upsilon).$$

**Proof.** This is a direct consequence of the decay of the minimizer established in Lemma 37. □

**Definition 41.** *Define*

$$\bar{G}^{f,\xi,\exp}(\upsilon) = \check{G}^{f,\xi,\exp}(\upsilon) \text{ if } \upsilon \in \mathcal{Y}^{M(f)}(\mu^{0,\xi}) \text{ and } +\infty \text{ otherwise.}$$

The next lemma discusses $\Gamma$-convergence of the variational description as the infrared cutoff is removed.

**Lemma 42.** *Let $(\xi_n)_n$ be a sequence of spatial cutoffs such that $\xi_n(x) = 1$ for $|x| \leqslant n$, and*

$$\sup_{n \in \mathbb{N}} \|\xi_n\|_{W^{k,\infty}} < +\infty, \qquad \text{for any } k \in \mathbb{N}.$$

*Pick a subsequence $\xi_n$ (not relabeled) such that $\mu^n \to \mu$ weakly on $\mathfrak{S} \times \hat{\mathfrak{H}}_w$ where $\mu^n$ is the minimizer of $\check{F}^{0,\xi_n}$ in $\bar{\mathcal{X}}$. Then $\bar{G}^{f,\xi_n,\exp}$ $\Gamma$-converges to $\bar{G}^{f,\infty,\exp}$ on $\mathcal{Y}^M$ where*

$$\bar{G}^{f,\infty,\exp}(\upsilon) = \mathbb{E}_\upsilon\left[ f(W + Z + K) + \int \exp(\beta Z)(\exp(\beta K) - 1) \mathrm{d}M^\beta + \mathcal{E}(K) \right]$$

*if $\upsilon \in \mathcal{Y}^M(\mu)$ and $+\infty$ otherwise.*

**Proof.** We first prove the liminf inequality. Assume that $\upsilon^n \to \upsilon$ in $\mathcal{Y}^M$. By definition of $\xi_n$ we can restrict ourselves to $\upsilon \in \mathcal{Y}^M(\mu)$. From the definition of convergence in $\mathcal{Y}^M$ it is not hard to see that

$$\mathbb{E}_{\upsilon^n}[f(W + Z^n + K^n)] \to \mathbb{E}_\upsilon[f(W + Z + K)],$$

and

$$\mathbb{E}_\upsilon[\mathcal{E}(K)] \leqslant \liminf_{n \to \infty} \mathbb{E}_{\upsilon^n}[\mathcal{E}(K^n)].$$

by the Portmanteau theorem. Now we can invoke the Skhorohod's embedding theorem from [24] and pass to a probability space $\tilde{\mathbb{P}}$ and on it a sequence $(W^n, M^{\beta,n}, Z^n, K)$ convergent almost surely on $\mathfrak{S} \times \hat{\mathfrak{H}}_w \times \mathfrak{H}_w(\rho^{-N})$ to $(W, M^\beta, Z, K)$ such that $\mathrm{Law}(W^n, M^{\beta,n}, Z^n, K) = \upsilon^n$. Using the fact that

$$\sup_n \mathbb{E}_{\tilde{\mathbb{P}}}\left[\|Z^n - Z\|^2_{\mathfrak{H}_w(\rho^{1/2})} + \|K^n - K\|^2_{\mathfrak{H}_w(\rho^{-N})}\right] < \infty$$

and almost-sure convergence we can assume that

$$\mathbb{E}_{\tilde{\mathbb{P}}}\|M^{\beta,n} - M^\beta\|^p_{B^\delta_{p,p}} \to 0, \qquad \mathbb{E}_{\tilde{\mathbb{P}}}\left[\|Z^n - Z\|^q_{\mathfrak{H}_w(\rho^{1/2})} + \|K^n - K\|^q_{\mathfrak{H}_w(\rho^{-N})}\right] \to 0,$$

for any $\delta, p$ satisfying the assumptions of Proposition 28 and $q < 2$. Let $\varsigma_N \colon \mathbb{R} \times \mathbb{R}^2 \to \mathbb{R}_+$ be continuous such that
$$\varsigma_N(z, x) = z \text{ if } z, |x| \leqslant N,$$
$\varsigma_N(z, x) = 0$ if either $z \geqslant 2N$ or $|x| \geqslant 2N$ and moreover $\varsigma_N(z, y) \leqslant \varsigma_{N+1}(z, y)$. We will write in the following $\varsigma_N(f, x)$ for the function $x \in \mathbb{R}^2 \mapsto \varsigma_N(f(x), x)$. We claim that
$$\lim_{n \to \infty} \int \xi_n \exp(\beta Z^n) \varsigma_N((\exp(\beta K^n) - 1)_+, x) \mathrm{d} M^{\beta, n}$$
$$= \int \exp(\beta Z) \varsigma_N((\exp(\beta K) - 1)_+, x) \mathrm{d} M^\beta.$$
Since
$$\sup_{n \in \mathbb{N}} \|\exp(\beta Z^n)\|_{L^\infty} \leqslant 1, \qquad \sup_{n \in \mathbb{N}} \|\varsigma_N((\exp(\beta K^n) - 1)_+, x)\|_{L^\infty} \leqslant N,$$
we can apply Lemma 45 below provided that, for any , $q < 2$,
$$\mathbb{E}\big[\|\xi_n \exp(\beta Z^n) \varsigma_N((\exp(\beta K^n) - 1)_+, x) - \exp(\beta Z) \varsigma_N((\exp(\beta K) - 1)_+, x)\|_{W^{1, q}(\rho^{-1})}^q\big] \to 0,$$
which in turn is ensured by Lemma 44 and Lemma 46. Then, recalling that $Z \leqslant 0$ $\mu^n$ almost surely and monotone convergence we have,
$$\int \exp(\beta Z)(\exp(\beta K) - 1)_+ \mathrm{d} M^\beta$$
$$= \sup_N \int \exp(\beta Z) \varsigma_N((\exp(\beta K) - 1)_+, x) \mathrm{d} M^\beta$$
$$= \sup_N \lim_{n \to \infty} \int \xi_n \exp(\beta Z^n) \varsigma_N((\exp(\beta K^n) - 1)_+, x) \mathrm{d} M^{\beta, n}$$
$$\leqslant \sup_N \liminf_{n \to \infty} \int \xi_n \exp(\beta Z^n)(\exp(\beta K^n) - 1)_+ \mathrm{d} M^{\beta, n}$$
$$= \liminf_{n \to \infty} \int \xi_n \exp(\beta Z^n)(\exp(\beta K^n) - 1)_+ \mathrm{d} M^{\beta, n}.$$
Let us now show that
$$\lim_{n \to \infty} \int \xi_n \exp(\beta Z^n)(\exp(\beta K^n) - 1)_- \mathrm{d} M^{\beta, n} = \int \exp(\beta Z)(\exp(\beta K) - 1)_- \mathrm{d} M^\beta.$$
To start, note that
$$|\exp(\beta Z^n)| \leqslant 1, \qquad |(\exp(\beta K^n) - 1)_-| \leqslant 1,$$
so again we have that from Lemmas 44, 46 below, for any $p < 2$,
$$\mathbb{E}\|\exp(\beta Z^n)(\exp(\beta K^n) - 1)_- - \exp(\beta Z)(\exp(\beta K) - 1)_-\|_{W^{1, p}(\rho^{-1})}^p \to 0.$$
From this and using that $\xi_n = 1$ for $|x| \leqslant n$, it is easy to see that
$$\mathbb{E}\|\xi_n \exp(\beta Z^n)(\exp(\beta K^n) - 1)_- - \exp(\beta Z)(\exp(\beta K) - 1)_-\|_{W^{1, p}(\rho^{-1/2})}^p \to 0.$$
we have that,
$$\int \xi_n \exp(\beta Z^n)(\exp(\beta K^n) - 1)_- \mathrm{d} M^{\beta, n} \to \int \exp(\beta Z)(\exp(\beta K) - 1)_- \mathrm{d} M^\beta.$$
To construct the recovery sequence consider $\upsilon$ such that $\bar{G}^{f, \infty, \exp}(\upsilon) < \infty$. Then
$$\mathbb{E}_\upsilon \int \exp(\beta Z)(\exp(\beta K) - 1) \mathrm{d} M^\beta < \infty, \qquad \mathbb{E}_\upsilon \mathcal{E}(K) \leqslant M$$

where the second inequality holds by definition of $\mathcal{Y}^M$. Again, after an application of Lemma 44 and interpolation between $L^\infty$ and $H^1(\rho^{-1})$ and provided $N > 1/\delta$, we have that

$$\mathbb{E}_v \int \exp(\beta Z)(\exp(\beta K) - 1)_- \mathrm{d}M^\beta$$
$$\leqslant \mathbb{E}\Big[\|(\exp(\beta K) - 1)_-\|_{B^\delta_{p',p'}(\rho^{-1})}^{p'\delta}\Big]\mathbb{E}\Big[\|M^\beta\|_{B^{-\delta}_{p,p}(\rho)}^p\Big]$$
$$\leqslant \mathbb{E}[\|(\exp(\beta K) - 1)_-\|_{H^1(\rho^{-N})}]\mathbb{E}\Big[\|M^\beta\|_{B^{-\delta}_{p,p}(\rho)}^p\Big]$$
$$< \infty.$$

So

$$\mathbb{E}_v \int \exp(\beta Z)(\exp(\beta K) - 1)_+ \mathrm{d}M^\beta$$
$$\leqslant \left|\mathbb{E}_v \int \exp(\beta Z)(\exp(\beta K) - 1)\mathrm{d}M^\beta\right| + \left|\mathbb{E}_v \int \exp(\beta Z)(\exp(\beta K) - 1)_-\mathrm{d}M^\beta\right|$$
$$< \infty.$$

Now set $K^N(x) = \varsigma_N(K, x)$, by monotone convergence

$$\lim_{N \to \infty} \mathbb{E}_v \int \exp(\beta Z)(\exp(\beta K^N) - 1)_+ \mathrm{d}M^\beta \to \mathbb{E}_v \int \exp(\beta Z)(\exp(\beta K) - 1)_+ \mathrm{d}M^\beta.$$

Moreover, since $|(\exp(\beta K^N) - 1)_-| \leqslant |(\exp(\beta K) - 1)_-|$, by dominated convergence also holds

$$\lim_{N \to \infty} \mathbb{E}_v \int \exp(\beta Z)(\exp(\beta K^N) - 1)_- \mathrm{d}M^\beta \to \mathbb{E}_v \int \exp(\beta Z)(\exp(\beta K) - 1)_- \mathrm{d}M^\beta.$$

Invoking [24] there exist a probability space $(\tilde{\Omega}, \tilde{\mathcal{F}}, \tilde{\mathbb{P}})$ and random variables $(W^n, M^\beta, Z^n)$ on it such that $(W^n, M^\beta, Z^n) \to (W, M^\beta, Z)$ almost surely in $\mathfrak{S} \times \mathfrak{H}_w(\rho)$ and $\mathrm{Law}_{\tilde{\mathbb{P}}}(W^n, M^\beta, Z^n) = \mu^n$.

Let $K$ be a random variable $\tilde{\Omega} \to \mathfrak{H}_w(\rho^{-1})$ such that $v = \mathrm{Law}_{\tilde{\mathbb{P}}}(W, Z, K)$. Set again $K^N = \varsigma^N(K)$ and $v^{n,N} = \mathrm{Law}_{\tilde{\mathbb{P}}}(W^n, Z^n, K^N)$. We claim that

$$\bar{G}^{f,\xi_n,\exp}(v^{n,N}) \to \bar{G}^{f,\infty,\exp}(v^N).$$

From Lemmas 44, 45, 46 we get

$$\mathbb{E}_{\tilde{\mathbb{P}}} \int \xi_n \exp(\beta Z^n)(\exp(\beta K^N) - 1)_- \mathrm{d}M^{\beta,n} \to \mathbb{E}_{\tilde{\mathbb{P}}} \int \exp(\beta Z)(\exp(\beta K^N) - 1)_- \mathrm{d}M^\beta.$$

Furthermore

$$(\exp(\beta K^N) - 1)_+ \leqslant \exp(N)\mathbb{1}_{B(0,2N)}, \qquad \exp(\beta Z^n) \leqslant 1,$$

and $(K^N, Z^n) \to (K^N, Z)$ in $H^1(B(0,N))$. So again from Lemmas 44, 45, 46

$$\mathbb{E}_{\tilde{\mathbb{P}}} \int \exp(\beta Z^n)(\exp(\beta K^N) - 1)_+ \mathrm{d}M^{\beta,n} \to \mathbb{E}_{\tilde{\mathbb{P}}} \int \exp(\beta Z)(\exp(\beta K^N) - 1)_+ \mathrm{d}M^\beta.$$

Finally

$$\lim_{n \to \infty} \mathbb{E}_{\tilde{\mathbb{P}}} f(W^n + Z^n + K^N) \to \mathbb{E}_{\tilde{\mathbb{P}}} f(W + Z + K^N)$$

follows by boundedness and continuity of $f$. Taking a diagonal sequence we obtain the statement. $\square$

**Corollary 43.** *If $(\xi_n)_n$ is as in Lemma 42 we have*

$$\lim_{n \to \infty} -\log \int \exp(-f(\varphi))\mathrm{d}\nu^{\xi_n,\exp} = \inf_{v \in \mathcal{Y}^M(\mu)} \bar{G}^{f,\infty,\exp}(v)$$

**Proof.** This follows from Γ-convergence, the fact that $\mathcal{Y}^M$ is sequentially compact and the fundamental theorem of Γ-convergence (see e.g. [14]). □

We conclude with some technical results needed in the proofs above.

**Lemma 44.** *Assume that $g: \mathbb{R} \to \mathbb{R}$ is continuously differentiable, $g, g'$ are bounded. Define $f = g_-$ (the negative part of $g$). We have that the map*

$$\varphi \mapsto f(\varphi),$$

*continuously sends $W_w^{1,2}(\rho^\alpha)$ in itself with*

$$\|f(\varphi)\|_{W^{1,2}(\rho^\alpha)} \lesssim 1 + \|\varphi\|_{W^{1,2}(\rho^\alpha)}.$$

*If $g(0) = 0$ we have in particular*

$$\|f(\varphi)\|_{W^{1,2}(\rho^\alpha)} \lesssim \|\varphi\|_{W^{1,2}(\rho^\alpha)}.$$

**Proof.** Assume that $g(0) = 0$. By definition of $f$ we have

$$\int |f(\varphi(x))|^2 \rho^\alpha(x) \mathrm{d}x \lesssim \int |\varphi(x)|^2 \rho^\alpha(x) \mathrm{d}x.$$

Furthermore by [15] p. 291–292, we have

$$\int |\nabla_x f(\varphi(x))|^2 \rho^\alpha(x) \mathrm{d}x = \int \mathbb{1}_{g(\varphi(x)) \leqslant 0} |g'(\varphi(x)) \nabla_x \varphi(x)|^2 \rho^\alpha(x) \mathrm{d}x \lesssim \int |\nabla_x \varphi(x)|^2 \rho^\alpha(x) \mathrm{d}x.$$

Now for weak continuity, let $\varphi_n \to \varphi$ weakly in $W^{1,2}(\rho^\alpha)$. We know that

$$\sup_n \|\varphi\|_{W^{1,2}(\rho^\alpha)} < \infty \text{ which implies } \sup_n \|f(\varphi)\|_{W^{1,2}(\rho^\alpha)} < \infty.$$

Now let $\varphi_{n_k}$ be a subsequence such that $f(\varphi_{n_k}) \to \xi$. We can extract a further subsequence (not relabeled) such that

$$\varphi_{n_k} \to \varphi$$

almost surely. Then $f(\varphi_{n_k}) \to f(\varphi)$ almost surely which in turn implies $\xi = f(\varphi)$. From this we can conclude. The second statement follows from the first by writing $g = \tilde{g} + g(0)$ with $\tilde{g}(0) = 0$. □

**Lemma 45.** *Let $M^{\beta,n}$ be a sequence of random measures in $\mathcal{P}(\mathbb{R}^2)$ on a probability space $(\Omega, \mathcal{F}, \mathbb{P})$, such that $M^{\beta,n} \to M^\beta$ $\mathbb{P}$-almost surely in $B_{p,p}^{-s}(\rho)$ for some $s > \frac{\beta^2}{4\pi}(p-1)$ and $p < 8\pi/\beta^2$ close enough to $\frac{\beta^2}{4\pi}(p-1), 8\pi/\beta^2$ respectively. There exists $q < 2$ such that we have*

$$\mathbb{E} \int f_n \mathrm{d}M^{\beta,n} \to \mathbb{E} \int f \mathrm{d}M^\beta$$

*whenever*

$$\mathbb{E}\bigl[\|f_n - f\|_{W^{1,q}(\rho^{-N})}^q\bigr] \to 0$$

*and*

$$\sup_{n \in \mathbb{N}} \|f_n\|_{L^\infty} + \|f\|_{L^\infty} = C < \infty.$$

**Proof.** By interpolation we have that for $\delta < q/p'$ and $N > \delta^{-1}$ we have with

$$\|f_n - f\|_{B_{p',p'}^\delta(\rho^{-1})} \leqslant C \|f_n - f\|_{W^{1,q}(\rho^{-N})}^\delta.$$

Now recall that for any $\frac{8\pi}{\beta^2} > p > 1$ and $\delta > \frac{\beta^2}{4\pi}(p-1)$ we have
$$\sup_n \mathbb{E}\Big[\|M^{\beta,n}\|^p_{B^{-\delta}_{p,p}(\rho)}\Big] < \infty.$$
so provided $1/p + 1/p' = 1$
$$\begin{aligned}
\mathbb{E}\Big[\int f_n \mathrm{d}M^{\beta,n} - \int f \mathrm{d}M^\beta\Big] &= \mathbb{E}\Big[\int (f_n - f)\mathrm{d}M^{\beta,n}\Big] + \mathbb{E}\Big[\int f \mathrm{d}M^{\beta,n} - \int f \mathrm{d}M^\beta\Big] \\
&\leqslant \mathbb{E}\Big[\|f_n - f\|_{B^\delta_{p',p'}(\rho^{-1})}\|M^{\beta,n}\|_{B^{-\delta}_{p,p}(\rho)}\Big] \\
&\quad + \mathbb{E}\Big[\|f\|_{B^\delta_{p',p'}(\rho^{-1})}\|M^{\beta,n} - M^\beta\|_{B^{-\delta}_{p,p}(\rho)}\Big] \\
&\leqslant \mathbb{E}\big[\|f_n - f\|^{p'\delta}_{W^{1,q}(\rho^{-1})}\big]^{\frac{1}{p'}} \mathbb{E}\Big[\|M^{\beta,n}\|^p_{B^{-\delta}_{p,p}(\rho)}\Big]^{\frac{1}{p}} \\
&\quad + \mathbb{E}\big[\|f\|^{p'\delta}_{W^{1,q}(\rho^{-1})}\big]^{\frac{1}{p'}} \mathbb{E}\Big[\|M^{\beta,n} - M^\beta\|^p_{B^{-\delta}_{p,p}(\rho)}\Big]^{\frac{1}{p}}. \\
&\to 0
\end{aligned}$$

since $p'\delta < q$. Note that for any $\beta^2 < 8\pi$ we can choose $q < 2$ such that
$$q/p' = \frac{q}{p}(p-1) > \frac{\beta^2}{4\pi}(p-1)$$
so if we choose $\delta$ such that
$$\frac{q}{p}(p-1) > \delta > \frac{\beta^2}{4\pi}(p-1)$$
we obtain as required
$$\frac{8\pi}{\beta^2} > p > 1 \qquad \delta > \frac{\beta^2}{4\pi}(p-1) \qquad \delta < q/p' = \frac{q}{p}(p-1). \qquad \square$$

**Lemma 46.** *For any $q < 2$ there exists an $N < \infty$ such that if*
$$\|f\|_{L^\infty}, \|g\|_{L^\infty} \leqslant C,$$
*then*
$$\|fg\|_{W^{1,p}(\rho^{-1})} \lesssim \|f\|_{H^1(\rho)}\|g\|_{H^1(\rho^{-N})}.$$
*Note that the $L^\infty$ norms are unweighted. If furthermore if $f_n \to f$ in $H^1(\rho)$ and $g_n \to g$ in $H^1(\rho^{-N})$ and $\sup_n \|f_n\|_{L^\infty} + \|g_n\|_{L^\infty} \leqslant C$ then*
$$f_n g_n \to fg \qquad \text{in } W^{1,q}(\rho^{-1}).$$

**Proof.** It is enough to prove the statement for smooth functions by density. By Hölder's inequality we have
$$\|fg\|_{L^p(\rho^{-1})} \leqslant \|f\|_{L^\infty}\|g\|_{L^2(\rho^{-N})}.$$
By the Leibniz rule
$$\|\nabla(fg)\|_{L^p(\rho^{-1})} \leqslant \|f\|_{L^\infty}\|\nabla g\|_{L^p(\rho^{-1})} + \|\nabla f\|_{L^2(\rho)}\|g\|_{L^p(\rho^{-\alpha N})},$$
provided that $\alpha N > 2$. By interpolation and provided that $\alpha/p = 1/2$, we also have
$$\|g\|_{L^p(\rho^{-\alpha N})} \leqslant \|g\|^{(1-\alpha)}_{L^\infty}\|g\|^\alpha_{L^2(\rho^{-N})} \leqslant C\|g\|^\alpha_{L^2(\rho^{-N})}.$$
Choosing $N$ large enough gives the statement. For the second statement we rewrite
$$fg - f_n g_n = g(f - f_n) + f_n(g - g_n)$$
and proceed analogously. $\square$

# Appendix A  Overview of the notation

| | |
|---:|:---|
| $\theta^0$ | the Gaussian Free Field with covariance $(m^2-\Delta)^{-1}$ |
| $\theta^V$ | Gibbs measure with potential $V$ |
| $\mathcal{W}^V(f)$ | $:=-\log\int\exp(-f(\phi))\,\theta^V(\mathrm{d}\phi)$ |
| $(W_s)_{s\in[0,1]}$ | Brownian martingale with covariance $(t\wedge s)(m^2-\Delta)^{-1}$ |
| $\rho$ | a spatial weight function |
| $B^s_{p,q}(\rho)$ | a weighted Besov space |
| $C^\alpha$ | Besov $-$ Hoelder space |
| $\mathcal{E}(Z)$ | $:=\dfrac{1}{2}\int_0^1\|(m^2-\Delta)^{1/2}\dot{Z}_s\|^2_{L^2}\mathrm{d}s$ (cost functional, entropy) |
| $G$ | a tame functional on $C([0,1],C^{-\varepsilon}(\rho^\varepsilon))$, see Definition 1. |
| $\mathfrak{H}(\rho^\alpha)$ | $=H^1_t H^1_x(\rho^\alpha)$ |
| $\mathfrak{H}$ | $=\mathfrak{H}(\rho^0)$ |
| $\mathfrak{H}^a(\rho^\alpha)$ | the space of processes in $\mathfrak{H}(\rho^\alpha)$ which are adapted to $W$ |
| $\mathfrak{H}_w(\rho^\alpha)$ | $\mathfrak{H}^a(\rho^\alpha)$ equipped with the weak topology |
| $\hat{\mathfrak{H}}$ | $:=\mathfrak{H}(\rho^{1/2})$ |
| $\mathfrak{S},\mathfrak{S}_{\exp}$ | space of enhanced noise with Wick powers GMC respectively |
| $\mu$ | measure on $\mathfrak{S}\times\mathfrak{H}^a(\rho^\alpha)$ |
| $\mathcal{X}$ | $:=\left\{\mu=\mathrm{Law}_{\mathbb{P}}(\mathbb{W},Z)\in\mathcal{P}(\mathfrak{S}\times\hat{\mathfrak{H}}_w):Z\in\hat{\mathfrak{H}}^a,\mathbb{E}_\mu[\|Z\|^2_{\hat{\mathfrak{H}}}]<\infty\right\}$ |
| $\bar{\mathcal{X}}$ | $:=\left\{\mu\in\mathcal{P}(\mathfrak{S}\times\hat{\mathfrak{H}}_w):\exists\mu_n\in\mathcal{X},\mu_n\to\mu\text{ weakly},\sup_{n\in\mathbb{N}}\mathbb{E}_{\mu_n}[\|Z\|^2_{\hat{\mathfrak{H}}}]<\infty\right\}$ |
| $F^G(Z)$ | $:=\mathbb{E}[G(W+Z)+\mathcal{E}(Z)]$ |
| $\check{F}^G(\mu)$ | $:=\mathbb{E}_\mu[G(W+Z)+\mathcal{E}(Z)]$ |
| $\bar{\mathcal{Y}}(\mu)$ | see Definition 4 |
| $\upsilon$ | element of $\mathcal{Y}(\mu)$ |
| $[\![\ ]\!]$ | Wick ordering |
| $\Lambda$ | compact subset of $\mathbb{R}^2$ |
| $V_\Lambda(\phi)$ | $:=\lambda\int_\Lambda[\![\phi^4]\!]$ |
| $\mu^{f+V}$ | minimizer of $F^{f+V}$ |
| $\beta\in(0,8\pi)$ | charge parameter |
| $M^\beta$ | the Gaussian multiplicative chaos (GMC) |
| $\xi\in C_c^\infty(\mathbb{R}^2)$ | a positive, compactly supported, spatial cutoff function |
| $V^\xi(\phi)$ | $=\lambda\int_{\mathbb{R}^2}\xi[\![\exp(\beta\phi)]\!]$ |